\documentclass[12pt]{article}

\DeclareMathAlphabet{\mathpzc}{OT1}{pzc}{m}{it}

\usepackage{amsmath}
\usepackage{amssymb}
\usepackage{amsfonts}
\usepackage{amsthm}
\usepackage{mathrsfs}
\usepackage{ifsym}
\usepackage{fontenc}
\usepackage{textcomp}
\usepackage{tipa}
\usepackage{enumerate}
\usepackage[pdftex]{color, graphicx}

\begin{document}

\title{{\bf  Solutions of diophantine equations as periodic points of $p$-adic algebraic functions, III}}       
\author{Patrick Morton}        
\date{April 21, 2020}          
\maketitle

\begin{abstract} All the periodic points of a certain algebraic function related to the Rogers-Ramanujan continued fraction $r(\tau)$ are determined.  They turn out to be $0, \frac{-1 \pm \sqrt{5}}{2}$, and the conjugates over $\mathbb{Q}$ of the values $r(w_d/5)$, where $w_d$ is one of a specific set of algebraic integers, divisible by the square of a prime divisor of 5, in the field $K_d=\mathbb{Q}(\sqrt{-d})$, as $-d$ ranges over all negative quadratic discriminants for which $\left(\frac{-d}{5}\right) = +1$.  This yields new insights on class numbers of orders in the fields $K_d$.  Conjecture 1 of Part I is proved for the prime $p=5$, showing that the ring class fields over fields of type $K_d$ whose conductors are relatively prime to $5$ coincide with the fields generated over $\mathbb{Q}$ by the periodic points (excluding -1) of a fixed $5$-adic algebraic function.
\end{abstract}

\section{Introduction.}

In Part I a periodic point of an algebraic function $w=\mathfrak{g}(z)$, with minimal polynomial $g(z,w)$ over $F(z)$, $F$ a given field (often algebraically closed), was defined to be an element $a$ of $F$, for which numbers $a_i \in F$ exist satisfying the simultaneous equations
$$g(a,a_1) = g(a_1,a_2) = \dots = g(a_{n-1},a) = 0,$$
for some $n \ge 1$.  The numbers $a_i = \mathfrak{g}(a_{i-1})$ in this definition are to be thought of as suitable values of the multi-valued function $\mathfrak{g}(z)$, determined by possibly different branches of $\mathfrak{g}(z)$ (when considered over $F=\mathbb{C}$).  Note that if the coefficients of $g(x,y)$ lie in a subfield $k$ of $F$, over which $F$ is algebraic, then the set of periodic points of $\mathfrak{g}(z)$ in $F$ is invariant under the action of $\textrm{Gal}(F/k)$.  In this part the main focus will be on the multi-valued function $\mathfrak{g}(z)$, whose minimal polynomial is the polynomial
$$g(x,y) = (y^4+2y^3+4y^2+3y+1)x^5-y(y^4-3y^3+4y^2-2y+1)$$
considered in Part II, related to the Rogers-Ramanujan continued fraction $r(\tau)$ (in the notation of \cite{du}).  Recall that the function $r(\tau)$ satisfies the modular equation
$$g(r(\tau),r(5\tau)) = 0, \ \ \tau \in \mathbb{H},$$
where $\mathbb{H}$ is the upper half-plane.  (See \cite{anb}, \cite{ber}, \cite{du}.) \medskip

I will show, that when transported to the $\mathfrak{p}$-adic domain -- specifically to $\textsf{K}_5(\sqrt{5})$, where $\textsf{K}_5$ is the maximal unramified algebraic extension of the $5$-adic field $\mathbb{Q}_5$ -- the ``multi-valued-ness" disappears, in that the $a_i = T_5^i(a) \in \textsf{K}_5(\sqrt{5})$ become values of a {\it single-valued} algebraic function $T_5(x)$, defined on a suitable domain $\textsf{D}_5 \subset \textsf{K}_5(\sqrt{5})$.  Thus, $5$-adically, $a$ and its companions $a_i$ are periodic points of $T_5(x)$ in the usual sense.  Setting $\varepsilon = \frac{-1+\sqrt{5}}{2}$, this single-valued algebraic function is given by the $5$-adically convergent series
\begin{equation}
T_5(x) = \ x^5 + 5 + \sqrt{5} \sum_{k=2}^\infty {a_k \left(\frac{5\sqrt{5}}{x^5-\varepsilon^5}\right)^{k-1}}, \ \ a_k =\sum_{j=1}^4{{j/5 \atopwithdelims ( ) k}},\\
\label{eq:1.1}
\end{equation}
for $x$ in the domain
$$\textsf{D}_5 = \{x \in \textsf{K}_5(\sqrt{5}): |x|_5 \le 1 \ \wedge \ x \not \equiv 2 \ (\textrm{mod} \ \sqrt{5})\}.$$
More precisely, half of the periodic points of $\mathfrak{g}(z)$ lie in $\textsf{D}_5$; namely, those which lie in the unramified extension $\textsf{K}_5$.  The other half are periodic points of the function $T \circ T_5^{-1} \circ T$ and lie in $T(\textsf{D}_5)$, where
$$T(x) = \frac{-(1+\sqrt{5})x+2}{2x+1+\sqrt{5}}.$$
The function $T_5(x)$ has the property that $y=T_5(x)$ is the unique solution in $\textsf{K}_5(\sqrt{5})$ of the equation $g(x,y)=0$, for any $x \in \textsf{K}_5(\sqrt{5})$ for which $x \not \equiv 2$ (mod $\sqrt{5}$).  Thus, $T_5(x)$ is one of the values of $\mathfrak{g}(x)$, for $x \in \textsf{D}_5$.  \medskip

In Part II \cite{mor8} it was shown that the conjugates over $\mathbb{Q}$ of the values $\eta = r(w/5)$ of the Rogers-Ramanujan continued fraction are periodic points of the algebraic function $\mathfrak{g}(z)$, for specific elements $w$ in the imaginary quadratic field $K=\mathbb{Q}(\sqrt{-d})$.  In this part it will be shown that these values are, together with $0$ and $\frac{-1\pm \sqrt{5}}{2}$, the {\it only} periodic points of $\mathfrak{g}(z)$.  Let $d_K$ denote the discriminant of $K=\mathbb{Q}(\sqrt{-d})$, where $\left(\frac{-d}{5}\right) = +1$, and let $\wp_5$ denote a prime divisor of $(5)=\wp_5 \wp_5'$ in $K$.  Recall that $p_d(x)$ is the minimal polynomial over $\mathbb{Q}$ of the value $r(w_d/5)$, where $w_d$ is given by equation (2) below. \bigskip

\noindent {\bf Theorem 1.1.} {\it The set of periodic points in $\overline{\mathbb{Q}}$ (or $\overline{\mathbb{Q}}_5$ or $\mathbb{C}$) of the multi-valued algebraic function $\mathfrak{g}(z)$ defined by the equation $g(z,\mathfrak{g}(z))=0$ consists of $0, \frac{-1\pm \sqrt{5}}{2}$, and the roots of the polynomials $p_d(x)$, for negative quadratic discriminants $-d=d_Kf^2$ satisfying $\left(\frac{-d}{5}\right)=+1$.  Over $\mathbb{C}$ the latter values coincide with the values $\eta = r(w_d/5)$ and their conjugates over $\mathbb{Q}$, where $r(\tau)$ is the Rogers-Ramanujan continued fraction and the argument $w_d \in K=\mathbb{Q}(\sqrt{-d})$ satisfies}
\begin{equation}
w_d=\frac{v+\sqrt{-d}}{2} \in R_K, \ \ \wp_5^2 \mid w_d, \ \ \textrm{and} \ (N(w_d),f) = 1.
\label{eq:1.2}
\end{equation}
{\it Over $\overline{\mathbb{Q}}_5$, all the periodic points of $\mathfrak{g}(z)$ lie in $\textsf{K}_5(\sqrt{5})$.  Moreover, the periodic points of $\mathfrak{g}(z)$ in $\textsf{K}_5$ are periodic points in $\textsf{D}_5$ of the single-valued $5$-adic function $T_5(x)$.}
\medskip

From this theorem and the results of Part II we can assert the following.  Let $F_d$ denote the abelian extension $F_d=\Sigma_5 \Omega_f$ ($d \neq 4f^2$) or $F_d=\Sigma_5 \Omega_{5f}$ ($d=4f^2 >4$) of $K=\mathbb{Q}(\sqrt{-d})$, where $\Sigma_5$ is the ray class field of conductor $(5)$ over $K$ and $\Omega_f$ is the ring class field of conductor $f$ over $K$.  Since $(f,5)=1$ and $\Omega_{5f} = \Omega_5 \Omega_f$ when $d \neq 4f^2$ (see \cite[Satz 3]{h}), then $F_d=\Sigma_5 \Omega_{5f}$ in either case.  Furthermore, $F_d$ coincides with what Cox \cite{co} calls the extended ring class field $L_{\mathcal{O},5}$ for the order $\mathcal{O}=\textsf{R}_{-d}$ of discriminant $-d$ in $K$.  Cox refers to Cho \cite{cho}, who denotes this field by $K_{(5),\mathcal{O}}$, but these fields are already discussed in S\"ohngen \cite[see p. 318]{so}, who shows they are generated by division values of the $\tau$-function, together with suitable values of the $j$-function.  See also Stevenhagen \cite{ste} and the monograph of Schertz \cite[p. 108]{sch}.  \bigskip

\noindent {\bf Theorem 1.2.} {\it Let $K=\mathbb{Q}(\sqrt{-d})$, with $\left(\frac{-d}{5}\right)=+1$ and $-d=d_K f^2$, as above.  If $\mathcal{O}=\textsf{R}_{-d}$ is the order of discriminant $-d$ in $K$, the extended ring class field $F_d = \Sigma_5 \Omega_{5f}$ over $K$ is generated over $\mathbb{Q}$ by a periodic point $\eta=r(w_d/5)$ of the function $\mathfrak{g}(z)$ ($w_d$ is as in (\ref{eq:1.2})), together with a primitive $5$-th root of unity $\zeta_5$:}
\begin{equation}
F_d = \Sigma_5\Omega_{5f} = \mathbb{Q}(\eta, \zeta_5).
\label{eq:1.3}
\end{equation}
{\it Conversely, if $\eta \neq 0, \frac{-1\pm \sqrt{5}}{2}$ is any periodic point of $\mathfrak{g}(z)$, then for some $-d=d_Kf^2$ for which $\left(\frac{-d}{5}\right)=+1$, the field $\mathbb{Q}(\eta, \zeta_5) = F_d$.  Furthermore, the field $\mathbb{Q}(\eta)$ generated by $\eta$ alone is the inertia field for the prime divisor $\wp_5$ or for its conjugate $\wp_5'$ in the field $F_d$.} \bigskip

This theorem provides explicit examples of Satz 22 in Hasse's {\it Zahlbericht} \cite{h1}, according to which any abelian extension of $K$ is obtained from $\Sigma = \Omega_f(\zeta_n)$, for some integer $f \ge 1$ and some $n$-th root of unity $\zeta_n$, by adjoining square-roots of elements of $\Sigma$.  This holds because $\eta = r(w_d/5)$ satisfies a quadratic equation over $\Omega_f(\zeta_5)$.  See \cite[Prop. 4.3, Cor. 4.7, Thm. 4.8]{mor8}.  \medskip

Here the method of Part I \cite{mor5} and \cite{mor7}, which yielded an interpretation and alternate derivation of special cases of a class number formula of Deuring, leads to the following {\it new} class number formula.  \bigskip

\noindent {\bf Theorem 1.3.}  {\it Let $\mathfrak{D}_{n,5}$ be the set of discriminants $-d=d_K f^2 \equiv \pm 1$ (mod $5$) of orders in imaginary quadratic fields $K=\mathbb{Q}(\sqrt{-d})$ for which the automorphism $\tau_5 = \left(\frac{F_{5,d}/K}{\wp_5}\right)$ has order $n$ in the Galois group $\textrm{Gal}(F_{5,d}/K)$, where $F_{5,d}$ is the inertia field for $\wp_5$ in the abelian extension $F_d/K$.  If $h(-d)$ is the class number of the order $\textsf{R}_{-d} \subset K$, then for $n>1$,}
\begin{equation}
\sum_{-d \in \mathfrak{D}_{n,5}}{h(-d)}=\frac{1}{2} \sum_{k \mid n}{\mu(n/k) 5^k}.
\label{eq:1.4}
\end{equation}

Based on this theorem and numerical calculations, I make the following \bigskip

\noindent {\bf Conjecture.}  {\it Let $q>5$ be a prime number.  Let  $L_{\mathcal{O},q}=L_{\textsf{R}_{-d},q}$ be the extended ring class field over $K=K_d=\mathbb{Q}(\sqrt{-d})$ for the order $\mathcal{O}=\textsf{R}_{-d}$ of discriminant $-d=d_K f^2$ in $K$, and let $h(-d)$ denote the class number of the order $\mathcal{O}$.  Then the following class number formula holds:
$$\sum_{-d \in \mathfrak{D}_{n,q}}{h(-d)}=\frac{2}{q-1} \sum_{k \mid n}{\mu(n/k) q^k}, \ \ n >1,$$
where $\mathfrak{D}_{n,q}$ is the set of discriminants $-d=d_K f^2$ for which $\left(\frac{-d}{q}\right)=+1$ and the Frobenius automorphism $\tau_q = \left(\frac{F_{q,d}/K_d}{\wp_q}\right)$ has order $n$, and $F_{q,d}$ is the inertia field for the prime divisor $\wp_q$ (dividing $q$ in $K_d$) in the abelian extension $L_{\textsf{R}_{-d},q}$ of $K_d$.} \bigskip

As was shown in \cite{mor8} for the prime $q=5$, the extension $L_{\textsf{R}_{-d},q}$ is equal to $\Sigma_q \Omega_f/K$, if $d \neq 3f^2$ or $4f^2$; and is equal to $\Sigma_q \Omega_{qf}/K$, if $q \equiv 1$ (mod $4$) and $d=4f^2$; or $q \equiv 1$ (mod $3$) and $d=3f^2$.  The field $F_{q,d}$ has degree $(q-1)/2$ and is cyclic over the ring class field $\Omega_f$ of conductor $f$ over $K$. \medskip

One naturally expects that this conjecture describes an aspect of a much more general phenomenon.  For example, one could consider families of quadratic fields $K=\mathbb{Q}(\sqrt{-d})$ for which the prime divisors $q$ of a given fixed integer $Q$ all split in $K$.  These are the $Q$-admissible quadratic fields.  Analogous formulas should hold for certain sets of class fields over the family of (imaginary?) abelian extensions of a fixed degree over $\mathbb{Q}$, whose Galois groups belong to a fixed isomorphism type, and in which a given rational prime $q$ splits. \medskip

In Section 6 I show that a similar situation exists for the algebraic function $w=\mathfrak{f}(z)$ whose minimal polynomial over $\overline{\mathbb{Q}}(z)$ is $h(z,w)$, where
\begin{align*}
h(z,w)=w^5&-(6+5z+5z^3+z^5)w^4+(21+5z+5z^3+z^5)w^3\\
&-(56+30z+30z^3+6z^5)w^2+(71+30z+30z^3+6z^5)w\\
&-120-55z-55z^3-11z^5.
\end{align*}
I showed in Part II (Theorem 5.4) that any ring class field $\Omega_f$ over the imaginary quadratic field $K$, whose conductor is relatively prime to $5$, is generated over $K$ by a periodic point $\upsilon$ of $\mathfrak{f}(z)$, which satisfies $\upsilon = \eta - \frac{1}{\eta}$, for a certain periodic point $\eta$ of $\mathfrak{g}(z)$.  In Theorem 6.2 of this paper I show that {\it any} periodic point $\upsilon \neq -1$ of $\mathfrak{f}(z)$ is related to a periodic point of $\mathfrak{g}(z)$ by $\upsilon = \eta - \frac{1}{\eta}=\phi(\eta)$, and that the $5$-adic function
$$\textsf{T}_5(x) = \phi \circ T_5 \circ \phi^{-1}(x), \ \ x \in \widetilde{\textsf{D}}_5=\phi(\textsf{D}_5 \cap  \{z \in \textsf{K}_5: |z|_5 = 1\}),$$
plays the same role for $\mathfrak{f}(z)$ that $T_5(x)$ plays for $\mathfrak{g}(z)$.  In particular, Theorems 6.2 and 6.3 show that Conjecture 1 of Part I is true for the prime $p=5$.  This leads to a proof of Deuring's formula for the prime $5$ in Theorem 6.5 and its corollary, analogous to the proof given in Part I and in \cite{mor7} for the prime $2$ and in \cite{mor4} for the prime $3$.

\section{Iterated resultants.}

Set
\begin{equation}
g(X,Y)=(Y^4+2Y^3+4Y^2+3Y+1)X^5-Y(Y^4-3Y^3+4Y^2-2Y+1).
\label{eq:2.1}
\end{equation}
In Part II \cite{mor8} it was shown that $(X,Y)=(\eta, \eta^{\tau_5})$, with $\eta=r(w_d/5)$ and $w_d$ given by (\ref{eq:1.2}), is a point on the curve $g(X,Y)=0$.  Here $\tau_5 = \left(\frac{\mathbb{Q}(\eta)/K}{\wp_5}\right)$ is the Frobenius automorphism for the prime divisor $\wp_5$ of $K=\mathbb{Q}(\sqrt{-d})$.  This fact implies that $r(w_d/5)$ and its conjugates over $\mathbb{Q}$ are periodic points of the function $\mathfrak{g}(z)$ defined by $g(z,\mathfrak{g}(z))=0$.  (See Part I, Theorem 5.3.)  In this section and Sections 3-4 it will be shown that these values, together with the fixed points $0, \frac{-1\pm \sqrt{5}}{2}$, represent {\it all} the periodic points of the algebraic function $\mathfrak{g}(z)$.  To do this we begin by considering a sequence of iterated resultants defined using the polynomial $g(x,y)$, as in Part I, Section 3. \medskip

We start by defining $R^{(1)}(x,x_1):=g(x,x_1)$, and note that
$$R^{(1)}(x,x_1) \equiv (x_1+3)^4(x^5-x_1) \ \ (\textrm{mod} \ 5).$$
Then we define the polynomial $R^{(n)}(x, x_n)$ inductively by
$$R^{(n)}(x,x_n):=\textrm{Resultant}_{x_{n-1}}(R^{(n-1)}(x,x_{n-1}),g(x_{n-1},x_n)), \ \ n \ge 2.$$
It is easily seen using induction that
$$R^{(n)}(x,x_n) \equiv (-1)^{n-1}(x_n+3)^{5^n-1}(x^{5^n}-x_n)  \ \ (\textrm{mod} \ 5),$$
so that the polynomial $R_n(x) :=R^{(n)}(x,x)$ satisfies
\begin{equation}
R_n(x) \equiv (-1)^{n-1}(x+3)^{5^n-1}(x^{5^n}-x)  \ \ (\textrm{mod} \ 5), \ \ n \ge 1.
\label{eq:2.3}
\end{equation}
The roots of $R_n(x)$ are all the periodic points of the multi-valued function $\mathfrak{g}(z)$ in any algebraically closed field containing $\mathbb{Q}$, whose periods are divisors of the integer $n$. (See Part I, p. 727.) \medskip

From this we deduce, by a similar argument as in the Lemma of Part I (pp. 727-728), that
$$\textrm{deg}(R_n(x))=2 \cdot 5^n-1, \ \ n \ge 1.$$
As in Part I, we define the expression $\textsf{P}_n(x)$ by
\begin{equation}
\textsf{P}_n(x) = \prod_{k \mid n}{R_k(x)^{\mu(n/k)}},
\label{eq:2.4}
\end{equation}
and show that $\textsf{P}_n(x) \in \mathbb{Z}[x]$.  From (\ref{eq:2.3}) it is clear that $R_n(x)$, for $n > 1$, is divisible (mod $5$) by the $N$ irreducible (monic) polynomials $\bar f_i(x)$ of degree $n$ over $\mathbb{F}_5$,
where
$$N= \frac{1}{n}\sum_{k \mid n}{\mu(n/k)5^k},$$
and that these polynomials are simple factors of $R_n(x)$ (mod $5$).  It follows from Hensel's Lemma that $R_n(x)$ is divisible by distinct irreducible polynomials $f_i(x)$ of degree $n$ over $\mathbb{Z}_5$, the ring of integers in $\mathbb{Q}_5$, for $1 \le i \le N$, with $f_i(x) \equiv \bar f_i(x)$ (mod $5$). In addition, all the roots of $f_i(x)$ are periodic of minimal period $n$ and lie in the unramified extension $\textsf{K}_5$.  Furthermore, $n$ is the smallest index $k$ for which $f_i(x) \mid R_k(x)$.  \medskip

Now we make use of the following identity for $g(x,y)$:
$$\left(x+\frac{1+\sqrt{5}}{2}\right)^5\left(y+\frac{1+\sqrt{5}}{2}\right)^5g(T(x),T(y)) = \left(\frac{5+\sqrt{5}}{2}\right)^5g(y,x),$$
where
$$T(x)=\frac{-(1+\sqrt{5})x+2}{2x+1+\sqrt{5}}.$$
We have
$$T(x)-2 = -\left(\frac{5+\sqrt{5}}{2}\right) \frac{2x-1+\sqrt{5}}{2x+1+\sqrt{5}}.$$
If the periodic point $a$ of $\mathfrak{g}(z)$, with minimal period $n > 1$, is a root of one of the polynomials $f_i(x)$, then $a$ is a unit in $\textsf{K}_5$, and for some $a_1, \dots, a_{n-1}$ we have
\begin{equation}
g(a,a_1)=g(a_1,a_2) = \cdots = g(a_{n-1},a)=0.
\label{eq:2.4a}
\end{equation}
Furthermore $a \not \equiv 2$ (mod  $\sqrt{5}$), since otherwise $a \equiv 2$ (mod $5$) would have degree $1$ over $\mathbb{F}_5$ (using that $\textsf{K}_5$ is unramified over $\mathbb{Q}_5$).  Hence, $2a+1+\sqrt{5}$ is a unit and $b=T(a) \equiv 2$ (mod $\sqrt{5}$).  All the $a_i$ satisfy $a_i \not \equiv 2$ (mod $\sqrt{5}$), as well, since the congruence $g(2,y) \equiv 4(y+3)^5$ (mod $5$) has only $y \equiv 2$ as a solution.  Hence, if some $a_i \equiv 2$, then $a_j \equiv 2$ for $j >i$, which would imply that $a \equiv 2$, as well.  The elements $b_i=T(a_i)$ are distinct and lie in $\textsf{K}_5(\sqrt{5})$, and the above identity implies that
\begin{equation}
g(b,b_{n-1})=g(b_{n-1},b_{n-2})= \dots = g(b_1,b)=0
\label{eq:2.5}
\end{equation}
in $\textsf{K}_5(\sqrt{5})$.  Thus, all the $b_i \equiv 2$ (mod $\sqrt{5}$), and the orbit $\{b,b_{n-1},\dots,b_1\}$ is distinct from all the orbits in (\ref{eq:2.4a}).  Now the map $T(x)$ has order $2$, so it is clear that $b=T(a)$ has minimal period $n$ in (\ref{eq:2.5}), since otherwise $a=T(b)$ would have period smaller than $n$.  It follows that there are at least $2N$ periodic orbits of minimal period $n>1$.  Noting that
$$R_1(x)= g(x,x) = \ x(x^2 + 1)(x^2 + x - 1)(x^4 + x^3 + 3x^2 - x + 1),$$
these distinct orbits and factors account for at least
\begin{align*}
2 \cdot 5-1+ \sum_{d \mid n,d>1}({2 \sum_{k \mid d}{\mu(d/k)5^k})} = -1+ 2\sum_{d \mid n}({\sum_{k \mid d}{\mu(d/k)5^k})} = 2 \cdot 5^n-1
\end{align*}
roots, and therefore all the roots, of $R_n(x)$.  This shows that the roots of $R_n(x)$ are distinct and the expressions $\textsf{P}_n(x)$ are polynomials.  Furthermore, over $\textsf{K}_5(\sqrt{5})$ we have the factorization
\begin{equation}
\textsf{P}_n(x) = \pm \prod_{1 \le i \le N}{f_i(x) \tilde f_i(x)}, \ \ n>1,
\label{eq:2.6}
\end{equation}
where $\tilde f_i(x) = c_i(2x+1+\sqrt{5})^{deg(f_i)}f_i(T(x))$, and the constant $c_i$ is chosen to make $\tilde f_i(x)$ monic.  Finally, the periodic points of $\mathfrak{g}(z)$ of minimal period $n$ are the roots of $\textsf{P}_n(x)$ and
\begin{equation}
\textrm{deg}(\textsf{P}_n(x)) = 2\sum_{k \mid n}{\mu(n/k)5^k}, \ \ n>1.
\label{eq:2.7}
\end{equation}
This discussion proves the following. \bigskip

\noindent {\bf Theorem 2.1.} {\it All the periodic points of $\mathfrak{g}(z)$ in $\overline{\mathbb{Q}}_5$ lie in $\textsf{K}_5(\sqrt{5})$.  The periodic points of minimal period $n$ coincide with the roots of the polynomial $\textsf{P}_n(x)$ defined by (\ref{eq:2.4}), and have degree $n$ over $\mathbb{Q}_5(\sqrt{5})$.  For $n>1$, exactly half of the periodic points of $\mathfrak{g}(z)$ of minimal period $n$ lie in $\textsf{K}_5$.} \bigskip

The last assertion in this theorem follows from the fact that $T(x)$ is a linear fractional expression in the quantity $\sqrt{5}$:
$$T(x) = \frac{-x\sqrt{5}-x+2}{\sqrt{5}+2x+1},$$
with determinant $-2(x^2+1)$.  If it were the case that $a \in \textsf{K}_5$ and $T(a) \in \textsf{K}_5$, for $n>1$, then the last fact would imply that $\sqrt{5} \in \textsf{K}_5$, which is not the case.  Therefore, for $n>1$, the only roots of $\textsf{P}_n(x)$ which lie in $\textsf{K}_5$ are the roots of the factors $f_i(x)$, in the above notation.  Furthermore, the factors $f_i(x)$ are irreducible over $\mathbb{Q}_5(\sqrt{5})$, since this field is purely ramified over $\mathbb{Q}_5$, which implies that the factors $\tilde f_i(x)$ are irreducible over $\mathbb{Q}_5(\sqrt{5})$, as well.

\section{A $5$-adic function.}

\noindent {\bf Lemma 3.1.} {\it Any root $\eta'$ of the polynomial $p_d(x)$ which is conjugate to $\eta = r(w_d/5)$ over $K = \mathbb{Q}(\sqrt{-d})$ satisfies $\eta' \not \equiv 2$} ({\it mod} $\mathfrak{p}$), {\it for any prime divisor $\mathfrak{p}$ of $\wp_5$} in $F_1 = \mathbb{Q}(\eta)$. \medskip

\noindent {\it Proof.}  It suffices to prove this for $\eta' = \eta$.  Assume $\eta \equiv 2$ (mod $\mathfrak{p}$), where $\mathfrak{p} \mid \wp_2$ in $F_1$.  Then the element $z=\eta^5-\frac{1}{\eta^5}$ satisfies $z \equiv 2^5-2^{-5} \equiv -1$ (mod $\mathfrak{p}$).  Hence the proof of \cite{mor8}, Theorem 4.6 implies that $d$ can only be one of the values $d=11,16,19$.  In these three cases $h(-d)=1$, so $\eta$ satisfies a quadratic polynomial over $K=\mathbb{Q}(\sqrt{-d})$.  We have
\begin{align*}
p_{11}(x) &= x^4-x^3+x^2+x+1\\
& = \left(x^2+\frac{-1+\sqrt{-11}}{2}x-1\right)\left(x^2+\frac{-1-\sqrt{-11}}{2}x-1\right);\\
p_{16}(x) &= x^4-2x^3+2x+1\\
& =(x^2+(-1+i)x-1)(x^2+(-1-i)x-1);\\
p_{19}(x) &= x^4+x^3+3x^2-x+1\\
& = \left(x^2+\frac{1+\sqrt{-19}}{2}x-1\right)\left(x^2+\frac{1-\sqrt{-19}}{2}x-1\right).
\end{align*}
In each case we have $\eta = r(w_d/5)$, where, respectively:
\begin{align*}
w_{11} &= \frac{33+\sqrt{-11}}{2}, \ N(w_{11}) = 5^2 \cdot 11,\\
w_{16} & = 3+4i, \ \ \ \ \ \ \ \ N(w_{16}) = 5^2,\\
w_{19} & = \frac{41+\sqrt{-19}}{2}, \ N(w_{19})=5^2 \cdot 17.
\end{align*}
Since $F_1=K(\eta)$ is unramified over $\wp_5$ and ramified over $\wp_5'$, the minimal polynomial $m_d(x)$ over $K$ of $\eta$ in each case is
the first factor listed above.  Since $\wp_5^2 \mid w_d$, we conclude that
$$\sqrt{-11} \equiv 2, \ i \equiv 3, \ \sqrt{-19} \equiv 4$$
modulo $\wp_5$ in $R_K$.  Then
$$m_{11}(x)  \equiv x^2+3x+4, \ m_{16}(x)  \equiv x^2+2x+4, \ m_{19}(x)  \equiv (x+1)(x+4)$$
modulo $\wp_5$, where the first two polynomials are irreducible mod $5$.  It follows that $\eta$ cannot be congruent to $2$ modulo any prime divisor of $\wp_5$.  In each case we also have $m_d(x) \equiv (x+3)^2$ (mod $\wp_5'$).  $\square$ \bigskip

Computing the partial derivative
\begin{align*}
\frac{\partial g(x,y)}{\partial y} = & \ (4y^3 + 6y^2 + 8y + 3)x^5 - 5y^4 + 12y^3 - 12y^2 + 4y - 1\\
\equiv & \ 4(x + 3)^5(y + 3)^3 \ (\textrm{mod} \ 5),
\end{align*}
we see that the points $(x,y)=(\eta,\eta^{\tau_5})$ on the curve $g(x,y)=0$ satisfy the condition
$$ \frac{\partial g(x,y)}{\partial y}|_{(x,y)=(\eta,\eta^{\tau_5})} \not \equiv 0 \ \textrm{mod} \ \mathfrak{p},$$
for any prime divisor $\mathfrak{p}$ of $\wp_5$.  Hence, the $\mathfrak{p}$-adic implicit function theorem implies that $\eta^{\tau_5}$ can be written as a single-valued function of $\eta$ in a suitable neighborhood of $x = \eta$.  (See \cite{pz}, p. 334.)  We shall now derive an explicit expression for this single-valued function.  \medskip

To do this, we consider $g(X,Y)=0$ as a quintic equation in $Y$.  Using Watson's method of solving a quintic equation from the paper \cite{lsw} of Lavallee, Spearman and Williams, we find that the roots $Y$ of $g(X,Y)=0$ are
\begin{align*}
Y=&\frac{Z+3}{5}+ \frac{\zeta}{10} (2Z+11+5\sqrt{5})^{4/5}(2Z+11-5\sqrt{5})^{1/5}\\
&+ \frac{\zeta^2}{10} (2Z+11+5\sqrt{5})^{3/5}(2Z+11-5\sqrt{5})^{2/5}\\
&+ \frac{\zeta^3}{10} (2Z+11+5\sqrt{5})^{2/5}(2Z+11-5\sqrt{5})^{3/5}\\
&+ \frac{\zeta^4}{10} (2Z+11+5\sqrt{5})^{1/5}(2Z+11-5\sqrt{5})^{4/5},
\end{align*}
where $\zeta$ is any fifth root of unity and $Z=X^5$.  This can also be written in the form
\begin{align*}
Y=&\frac{Z+3}{5}+ \frac{\zeta}{5} (Z-\bar \varepsilon^5)^{4/5}(Z-\varepsilon^5)^{1/5} + \frac{\zeta^2}{5} (Z-\bar \varepsilon^5)^{3/5}(Z- \varepsilon^5)^{2/5}\\
& + \frac{\zeta^3}{5} (Z-\bar \varepsilon^5)^{2/5}(Z- \varepsilon^5)^{3/5}+ \frac{\zeta^4}{5} (Z-\bar \varepsilon^5)^{1/5}(Z- \varepsilon^5)^{4/5},\\
=& \frac{Z+3}{5}+\frac{1}{5}(Z-\varepsilon^5)(U^4+U^3+U^2+U), \ \ U =\zeta^{-1} \left(\frac{Z-\bar \varepsilon^5}{Z-\varepsilon^5}\right)^{1/5}.
\end{align*}

Now, $\varepsilon^5 =\frac{-11+5\sqrt{5}}{2} \equiv \frac{-1}{2} \equiv 2$ (mod $5$), so for $\zeta=1$ and $Z \not \equiv 2$ (mod $5$), the functions $U^j$ can be expanded into a convergent series:
$$U^j=\left(\frac{Z-\bar \varepsilon^5}{Z-\varepsilon^5}\right)^{j/5}=\left(1+\frac{\varepsilon^5-\bar \varepsilon^5}{Z-\varepsilon^5}\right)^{j/5}= \sum_{k=0}^\infty{{\frac{j}{5} \atopwithdelims ( ) k}\left(\frac{5\sqrt{5}}{Z-\varepsilon^5}\right)^k}.$$
This series converges for all $Z \not \equiv 2$ (mod $\sqrt{5}$) in the field $\textsf{K}_5(\sqrt{5})$.  The terms in this series tend to $0$ in the $5$-adic valuation, because
$$5^k {\frac{j}{5} \atopwithdelims ( ) k} = \frac{j(j-5)(j-10) \cdots (j-5(k-1))}{k!}$$
and because the additive $5$-adic valuation of $k!$ satisfies
$$v_5(k!) = \frac{k-s_k}{4} \le \frac{k}{4},$$
where $s_k$ is the sum of the $5$-adic digits of $k$.  Thus, for all $x \not \equiv 2$ (mod $\sqrt{5}$) in $\textsf{K}_5(\sqrt{5})$ the expression
\begin{equation}
y=T_5(x)=\frac{x^5+3}{5}+\frac{1}{5}(x^5-\varepsilon^5)\sum_{k=0}^\infty {a_k \left(\frac{5\sqrt{5}}{x^5-\varepsilon^5}\right)^k}, \ \ a_k =\sum_{j=1}^4{{\frac{j}{5} \atopwithdelims ( ) k}},
\label{eq:3.1}
\end{equation}
represents a root of the equation $g(x,y)=0$ in the field $\textsf{K}_5(\sqrt{5})$.  This formula for $T_5(x)$ simplifies to:
\begin{equation}
T_5(x) = \ x^5 + 5 + \sqrt{5} \sum_{k=2}^\infty {a_k \left(\frac{5\sqrt{5}}{x^5-\varepsilon^5}\right)^{k-1}}.\\
\label{eq:3.2}
\end{equation}
Note that
\begin{equation}
T_5(x) \equiv x^5 \ (\textrm{mod} \ 5), \ \ |x|_5 \le 1.
\label{eq:3.3}
\end{equation}
This follows from the fact that $5$ divides the individual terms $b_k=5^k a_k (\sqrt{5})^{k-2}$ (ignoring the unit denominators) in the series (\ref{eq:3.2}), for $2 \le k \le 7$, as can be checked by direct computation, and from the following estimate for $v_5(b_k)$, the normalized additive valuation of $b_k$ in $\textsf{K}_5(\sqrt{5})$:
$$v_5(5^k a_k (\sqrt{5})^{k-2}) \ge \frac{k}{2}-1-\frac{k}{4} = \frac{k}{4}-1 \ge 1, \ \textrm{for} \ k \ge 8.$$
It follows from this that the function $T_5(x)$ can be iterated on the set
\begin{equation}
\textsf{D}_5=\{x \in \textsf{K}_5(\sqrt{5}): |x|_5 \le 1 \ \wedge \ x \not \equiv 2 \ (\textrm{mod} \ \sqrt{5})\}.
\label{eq:3.4}
\end{equation}

I claim now that (\ref{eq:3.1}) (or (\ref{eq:3.2})) gives the {\it only} root of $g(x,y)=0$ in the field $\textsf{K}_5(\sqrt{5})$, for a fixed $x \not \equiv 2$ (mod $\sqrt{5}$).  From the above formulas, a second root of this equation must have the form
$$y_1=\frac{x^5+3}{5}+\frac{1}{5}(x^5-\varepsilon^5)(U^4+U^3+U^2+U),$$
where
$$U =\zeta^{-1} \left(\frac{x^5-\bar \varepsilon^5}{x^5-\varepsilon^5}\right)^{1/5},$$
for some fifth root of unity $\zeta$.  But then 
$$U^4+U^3+U^2+U = \frac{U^5-1}{U-1}-1 \in \textsf{K}_5(\sqrt{5}),$$
so $U \in \textsf{K}_5(\sqrt{5})$; and since $\zeta U$ is also in $\textsf{K}_5(\sqrt{5})$, it follows that $\zeta \in \textsf{K}_5(\sqrt{5})$.  This is impossible, since the ramification index of $5$ in $\textsf{K}_5(\zeta)$ is $e=4$, while the ramification index of $5$ in $\textsf{K}_5(\sqrt{5})$ is only $e=2$. \bigskip

\noindent {\bf Proposition 3.2.}
{\it If $x \in \textsf{D}_5$, the subset of $\textsf{K}_5(\sqrt{5})$ defined by} (\ref{eq:3.4}), {\it then the series
\begin{equation}
y=T_5(x)=\ x^5 + 5 + \sqrt{5} \sum_{k=2}^\infty {a_k \left(\frac{5\sqrt{5}}{x^5-\varepsilon^5}\right)^{k-1}}, \ \ a_k =\sum_{j=1}^4{{\frac{j}{5} \atopwithdelims ( ) k}},
\label{eq:3.5}
\end{equation}
gives the unique solution of the equation $g(x,y)=0$ in the field $\textsf{K}_5(\sqrt{5})$.  Moreover, the image $T_5(x)$ also lies in $\textsf{D}_5$, so the map $T_5$ can be iterated on this set.} \bigskip

\noindent {\bf Corollary 3.3.} {\it The function $T_5(x)$ satisfies $T_5(\textsf{D}_5 \cap \textsf{K}_5) \subseteq \textsf{D}_5 \cap \textsf{K}_5$.} \medskip

\noindent {\it Proof.} Let $\sigma$ denote the non-trivial automorphism of $\textsf{K}_5(\sqrt{5})/\textsf{K}_5$.  If $x \in \textsf{D}_5 \cap \textsf{K}_5$, then $g(x,T_5(x))=0$ and $T_5(x) \in \textsf{K}_5(\sqrt{5})$ imply that $g(x^\sigma,T_5(x)^\sigma)=g(x,T_5(x)^\sigma) = 0$.  The theorem gives that $T_5(x)^\sigma = T_5(x)$, implying that $T_5(x) \in \textsf{K}_5$.  $\square$ \bigskip

Now the completion $(F_1)_\mathfrak{p}$ of the field $F_1 = \mathbb{Q}(\eta)$ with respect to a prime divisor $\mathfrak{p}$ of $R_{F_1}$ dividing $\wp_5$ is a subfield of $\textsf{K}_5(\sqrt{5})$.  This is because $F_1$ is unramified at the prime $\mathfrak{p}$ and is abelian over $K$, so that $(F_1)_\mathfrak{p}$ is unramified and abelian over $K_{\wp_5} = \mathbb{Q}_5$.  \medskip

By Lemma 3.1, we can substitute $x = \eta$ in (\ref{eq:3.5}), and since $\eta^{\tau_5}$ is a solution of $g(\eta,Y)=0$ in $\textsf{K}_5$, we conclude that $\eta^{\tau_5}=T_5(\eta)$.  Letting $\zeta = 1$ and $U=-u$ gives
$$\eta^{\tau_5}=\frac{\eta^5+3}{5}+\frac{1}{5}(\eta^5-\varepsilon^5)(u^4-u^3+u^2-u), \ \ u = -\left(\frac{\eta^5-\bar \varepsilon^5}{\eta^5-\varepsilon^5}\right)^{1/5}=\frac{1}{\varepsilon \xi} \in F;$$
which agrees with the result of \cite{mor8}, Theorem 5.1.  The automorphism $\tau_5$ is canonically defined on the unramified extension $\mathbb{Q}_5(\eta)$; defining $\tau_5$ to be trivial on $\mathbb{Q}_5(\sqrt{5})$, we have that $T_5(\eta^{\tau_5}) = T_5(\eta)^{\tau_5}$, and hence that
\begin{equation}
\eta^{\tau_5^n}=T_5^n(\eta), \ \ n \ge 1.
\label{eq:3.6}
\end{equation}
This also follows inductively from
$$g(\eta^{\tau_5^{n-1}},\eta^{\tau_5^n}) = g(\eta^{\tau_5^{n-1}},T_5(\eta^{\tau_5^{n-1}})) = g(\eta^{\tau_5^{n-1}},T_5^n(\eta))= 0.$$  Therefore, $\eta = r(w/5)$ is a periodic point of $T_5$ in $\textsf{D}_5$, and the minimal period of $\eta$ with respect to $T_5$ is equal to the order of the automorphism $\tau_5=\left(\frac{F_1/\Omega_f}{\wp_5}\right)$. \medskip

By Theorem 2.1, the periodic points of $\mathfrak{g}(z)$ lie in $\textsf{K}_5(\sqrt{5})$.  In particular, the minimal period of $\eta=r(w_d/5)$ with respect to $\mathfrak{g}(z)$ is the order $n$ of the automorphism $\tau_5$.  This is because any  values $\eta_i$, for which
$$g(\eta,\eta_1) = g(\eta_1, \eta_2) = \dots = g(\eta_{m-1}, \eta) = 0,$$
must themselves be periodic points with $\eta_i \not \equiv 2$ (mod $\sqrt{5}$).  This implies that $\eta_i \in \textsf{D}_5$, and then $\eta_i = T_5^i(\eta)$ follows from Proposition 3.2, so that $m$ must be a multiple of $n$.  Hence, $\eta = r(w_d/5)$ must be a root of the polynomial $\textsf{P}_n(x)$.  \bigskip

\noindent {\bf Theorem 3.4.}  {\it For any discriminant $-d \equiv \pm 1$ (mod $5$), for which the automorphism $\tau_5 = \left(\frac{F_1/K}{\wp_5}\right)$ has order $n$, the polynomial $p_d(x)$ divides $\textsf{P}_n(x)$.}

\section{Identifying the factors of $P_n(x)$.}

We will now show that the polynomials $p_d(x)$ in Theorem 3.4 are the only irreducible factors of $\textsf{P}_n(x)$ over $\mathbb{Q}$.  The argument is similar to the argument in \cite[pp. 877-878]{mor4} with added complexity due to the nontrivial nature of the points in $E_5[5]-\langle (0,0) \rangle$, plus the necessity of dealing with the action of the icosahedral group in this case. \medskip

To motivate the calculation below, we prove the following lemma. As in Part II, $F_1$ denotes the field $F_1 = \mathbb{Q}(\eta)$, where $\eta = r(w_d/5)$. \bigskip

\noindent {\bf Lemma 4.1.} {\it If $w=w_d$ is defined as in (\ref{eq:1.2}), and $\tau_5 = \left(\frac{F_1/\mathbb{Q}}{\wp_5}\right)$, then for some $5$-th root of unity $\zeta^i$, we have}
$$\eta^{\tau_5^{-1}} = r\left(\frac{w}{5}\right)^{\tau_5^{-1}}= \zeta^i r\left(\frac{w}{25}\right).$$ \smallskip

\noindent {\it Proof.}  Define $\tau_5$ on $F_1(\sqrt{5})=\mathbb{Q}(\eta,\sqrt{5})$ so that it fixes $\sqrt{5}$.  This is possible since $F_1$ and $K(\sqrt{5})$ are disjoint, abelian extensions of $K$.  (See the discussion in Sections 5.2 and 5.3 of \cite{mor8}, where $\tau_5=\sigma_1\phi|_{F_1}$ and both $\sigma_1$ and $\phi$ fix the field $L=\mathbb{Q}(\zeta)$.)  Recall from Part II that
$$\tau(b) = \frac{-b+\varepsilon^5}{\varepsilon^5 b+1}.$$
From $\tau(\xi^5)= \eta^5$ and $T(\eta^{\tau_5}) = \xi$ (Part II, Thms. 3.3 and 5.1) we then obtain
$$\eta^{5\tau_5^{-1}}=\tau(\xi^5)^{\tau_5^{-1}} = \tau\left((\xi^{\tau_5^{-1}})^5\right) = \tau(T(\eta)^5)=\mathfrak{r}(\eta),$$
where
$$\mathfrak{r}(z)=z\frac{z^4-3z^3+4z^2-2z+1}{z^4+2z^3+4z^2+3z+1},$$
as in the Introduction to Part II.  On the other hand,
$$\mathfrak{r}(\eta) = \mathfrak{r} \left(r\left(\frac{w}{5}\right)\right) = r^5\left(\frac{w}{25}\right),$$
by Ramanujan's modular equation.  Thus, $\eta^{5\tau_5^{-1}}=r^5(w/25)$, and the assertion follows.  $\square$
\bigskip

By (\ref{eq:3.3}), we have $f_i(T_5(x)) \equiv f_i(x^5)$ (mod $5$), and since $T_5(a)$ is an ''unramified'' periodic point in $\textsf{D}_5$ whenever $a$ is, it follows that $\sigma: x \rightarrow T_5(x)$ is a lift of the Frobenius automorphism on the roots of $f_i(x)$, for each $i$ with $1 \le i \le N$.  We may assume that $\sigma$ fixes $\sqrt{5}$, since $\textsf{K}_5$ and $\mathbb{Q}_5(\sqrt{5})$ are linearly disjoint over $\mathbb{Q}_5$.  In order to apply $\sigma$ to all the maps occurring in the proof below, we also extend $\sigma$ to the field $\textsf{K}_5 \mathbb{Q}_5\left(\sqrt{\frac{-5+\sqrt{5}}{2}}\right)$, so that it fixes elements of the field $\mathbb{Q}_5\left(\sqrt{\frac{-5+\sqrt{5}}{2}}\right)$, which is a cyclic quartic and totally ramified extension of $\mathbb{Q}_5$ (the minimal polynomial of the square-root being the Eisenstein polynomial $x^4 + 5x^2 + 5$).
\bigskip

\noindent {\bf Theorem 4.2.} {\it For $n > 1$ the polynomial $\textsf{P}_n(x)$ is a product of polynomials $p_d(x)$:
\begin{equation}
\textsf{P}_n(x)=\pm \prod_{-d \in \mathfrak{D}_{n,5}}{p_d(x)},
\label{eq:4.1}
\end{equation}
where $\mathfrak{D}_{n,5} $ is the set of discriminants $-d=d_K f^2$ of imaginary quadratic orders $\textsf{R}_{-d} \subset K = \mathbb{Q}(\sqrt{-d})$ for which $\left(\frac{-d}{5}\right)=+1$ and the corresponding automorphism $\tau_5 = \left(\frac{F_1/K}{\wp_5}\right)$ has order $n$ in $\textrm{Gal}(F_1/K)$.  Here $F_1=\mathbb{Q}(r(w_d/5))$ is the inertia field for the prime divisor $\wp_5 = (5, w_d)$ in the abelian extension $\Sigma_5 \Omega_f$ ($d \neq 4f^2$) or $\Sigma_5 \Omega_{5f}$ ($d=4f^2 > 4$) of $K$; and $p_d(x)$ is the minimal polynomial of the value $r(w_d/5)$ over $\mathbb{Q}$.}
\medskip

\noindent {\it Proof.} Let $\eta=\eta_0,\eta_1, \dots, \eta_{n-1}$, $n \ge 2$, be a periodic orbit of $T_5(x)$ contained in $\textsf{D}_5$, where $T_5^n(\eta) = \eta$, and let
$$\xi=T(\eta_1)=T(T_5(\eta))=T(\eta^\sigma).$$
Then the relation $g(\eta,\eta_1) = g(\eta,T(\xi))=0$ implies that $(\eta,\xi)$ is a point on the curve
$$\mathcal{C}_5: \ X^5+Y^5 = \varepsilon^5(1-X^5 Y^5).$$
Rewrite this relation as
$$\xi^5 = \frac{-\eta^5+\varepsilon^5}{\varepsilon^5 \eta^5+1} = \tau(\eta^5), \ \ \tau(b) = \frac{-b+\varepsilon^5}{\varepsilon^5 b+1}.$$
Let
$$E_5(b): \ Y^2 + (1+b)XY+bY=X^3+bX^2$$
be the Tate normal form for a point of order $5$; and let $E_{5,5}(b)$ be the isogenous curve
\begin{align*}
E_{5,5}(b): \ Y^2+(1+b)XY+5bY=X^3 &+7bX^2+6(b^3+b^2-b)X\\
& + b^5+b^4-10b^3-29b^2-b.
\end{align*}
The $X$-coordinate of the map $\psi:E_5(b) \rightarrow E_{5,5}(b)$ is given by
$$X(\psi(P)) = \frac{b^4+(3b^3+b^4)x+(3b^2+b^3)x^2+(b-b^2-b^3)x^3+x^5}{x^2(x+b)^2}, \ \ b=\eta^5,$$
with $x=X(P)$.  Note that $\textrm{ker}(\psi)=\langle (0,0) \rangle$, and $\psi$ is defined over $\mathbb{Q}(b)$.  (See  \cite[p. 259]{mor}.)  \medskip

The relation $\xi^5=\tau(\eta^5)$ implies that there is an isogeny $\phi: E_5(\eta^5) \rightarrow E_5(\tau(\eta^5)) = E_5(\xi^5)$.  This is because the $j$-invariant of $E_5(\xi^5)$ is
\begin{align*}
j_\xi &= \ \frac{(1-12\xi^5+14\xi^{10}+12\xi^{15}+\xi^{20})^3}{\xi^{25}(1-11\xi^5-\xi^{10})}\\
& = \ \frac{(1+228\eta^5+494\eta^{10}-228\eta^{15}+\eta^{20})^3}{\eta^5(1-11\eta^5-\eta^{10})^5},
\end{align*}
where the latter value is $j(E_{5,5}(\eta^5))$.  Thus, $E_{5,5}(\eta^5) \cong E_5(\xi^5)$ by an isomorphism $\iota_1$.  Composing $\psi$ (for $b=\eta^5$) with this isomorphism gives the isogeny $\phi=\iota_1 \circ \psi$.  Furthermore, $j(E_{5,5}(\eta^5))$ is invariant under the substitution $\eta \rightarrow T(\eta)=\xi^{\sigma^{-1}}$, so
\begin{align*}
j_\xi & = \ \left(\frac{(1+228\xi^5+494\xi^{10}-228\xi^{15}+\xi^{20})^3}{\xi^5(1-11\xi^5-\xi^{10})^5}\right)^{\sigma^{-1}}\\
 & = \ \left(\frac{(1-12\eta^5+14\eta^{10}+12\eta^{15}+\eta^{20})^3}{\eta^{25}(1-11\eta^5-\eta^{10})}\right)^{\sigma^{-1}}\\
 & = j_{\eta^{\sigma^{-1}}}.
 \end{align*}
 It follows that $E_5(\xi^5) \cong E_5((\eta^{\sigma^{-1}})^5)$ by an isomorphism $\iota_2$.  Composing $\iota_2$ with $\phi$ gives an isogeny $\iota_2 \circ \phi = \phi_1: E_5(\eta^5) \rightarrow E_5(\eta^5)^{\sigma^{-1}}$ of degree $5$.  Applying $\sigma^{-i+1}$ to the coefficients of $\phi_1$ gives an isogeny
$$\phi_i: \  E_5(\eta^5)^{\sigma^{-(i-1)}} \rightarrow E_5(\eta^5)^{\sigma^{-i}}, \ \ 1 \le i \le n,$$
which also has degree $5$.  Hence, $\iota = \phi_n \circ \phi_{n-1} \circ \cdots \circ \phi_1$ is an isogeny from $E_5(\eta^5)$ to $E_5(\eta^5)^{\sigma^{-n}}$ of degree $5^n$.  But $\sigma^n$ is trivial on $\mathbb{Q}_5(\eta, \sqrt{5})$, since $T_5^n(\eta) = \eta$.  Hence, $\iota: E_5(\eta^5) \rightarrow E_5(\eta^5)$.  \medskip
 
We will show that $\iota$ is a cyclic isogeny by showing that some point $P \in E_5(\eta^5)[5]$ is not in $\text{ker}(\iota)$.  The following formula from \cite{mor6} gives the $X$-coordinate on $E_5(b)$ for a point $P$ of order $5$, which does not lie in $\langle (0,0) \rangle$:
\begin{align*}
X(P)=\frac{-\varepsilon^4}{2}\frac{(- 2u^2+(1+\sqrt{5})u - 3\sqrt{5} - 7)(2u^2 + (2\sqrt{5} + 4)u + 3\sqrt{5} + 7)}{(- 2u^2+(\sqrt{5}+1)u - 2)(u + 1)^2},
\end{align*}
where
$$u^5 = -\frac{b-\bar \varepsilon^5}{b-\varepsilon^5}, \ \ b = \eta^5, \ \bar \varepsilon =-\frac{1+\sqrt{5}}{2}.$$
A calculation on Maple shows that
$$X_1=X(\psi(P))= \frac{-5+\sqrt{5}}{10}(b^2+\varepsilon^4 b+\bar \varepsilon^2), \ \ b=\eta^5.$$
This is the $X$-coordinate of the point $P'=\psi(P)$ on $E_{5,5}(b)$.  On the other hand, an isomorphism $\iota_1: E_{5,5}(b) \rightarrow E_5(\tau(b))$ is given by $\iota_1(X_1,Y_1) = (X_2,Y_2)$, where
$$X_2= \lambda_1^2 X_1+\lambda_1^2 \frac{b^2+30b+1}{12}-\frac{\tau(b)^2+6\tau(b)+1}{12},$$
and
$$\lambda_1^2 = \frac{\sqrt{5}\bar \varepsilon^5}{(b-\bar \varepsilon^5)^2}=\frac{\sqrt{5}\bar \varepsilon^5}{(\eta^5-\bar \varepsilon^5)^2}.$$
Under this isomorphism, $X_1=X(\psi(P))$ maps to $X_2=0$, whence $\phi(P)=\iota_1 \circ \psi(P) = \pm(0,0)$ on $E_5(\tau(b))=E_5(\xi^5)$.  Note that the map $\phi$ is defined over $\Lambda=\mathbb{Q}\left(\eta,\sqrt{\sqrt{5} \bar \varepsilon}\right)=\mathbb{Q}\left(\eta, \sqrt{\frac{-5-\sqrt{5}}{2}}\right)$, since $\lambda_1$ lies in this field.
\medskip

Now we find an explicit formula for the isomorphism $\iota_2$ between $E_5(\xi^5)$ and $E_5(\eta^{5\sigma^{-1}})$.  The Weierstrass normal form $Y^2=4X^3-g_2X-g_3$ of $E_5(b)$ has coefficients
\begin{align*}
g_2(b)=& \ \frac{1}{12}(b^4+12b^3+14b^2-12b+1),\\
g_3(b) = & \ \frac{-1}{216}(b^2+1)(b^4+18b^3+74b^2-18b+1).
\end{align*}
An isomorphism $\iota_2: \ E_5(\xi^5) \rightarrow E_5(\eta^{5\sigma^{-1}})$ is determined by a number $\lambda_2$ satisfying the equations
$$g_2(\eta^{5\sigma^{-1}})=\lambda_2^4 \cdot g_2(\xi^5), \ \ g_3(\eta^{5\sigma^{-1}})=\lambda_2^6 \cdot g_3(\xi^5).$$
We now use computations analogous to those in Lemma 4.1, obtaining
$$\eta^{5\sigma^{-1}}=\tau(\xi^5)^{\sigma^{-1}} = \tau\left((\xi^{\sigma^{-1}})^5\right) = \tau(T(\eta)^5)=\mathfrak{r}(\eta).$$
Then we solve for $\lambda_2^2$ from
$$\lambda_2^2 = \frac{g_3(\mathfrak{r}(\eta)) g_2(\tau(\eta^5))}{g_2(\mathfrak{r}(\eta))g_3(\tau(\eta^5))}$$
and find that
$$\lambda_2^2 = \frac{(11\sqrt{5} - 25)(2\eta + 1 + \sqrt{5})^2(- 2\eta^2 +(3+ \sqrt{5})\eta - 3 - \sqrt{5})^2}{40(-2\eta^2 - 2\eta - 3 + \sqrt{5})^2}.$$
Here, $\lambda_2$ lies in the field $\mathbb{Q}\left(\eta, \sqrt{-\sqrt{5} \varepsilon}\right)=\mathbb{Q}\left(\eta, \sqrt{\frac{-5+\sqrt{5}}{2}}\right)$, which coincides with the field $\Lambda$ above.  Hence, the desired isomorphism is given on $X$-coordinates by
$$X_3=\iota_2(X_2)=\lambda_2^2 X_2+\lambda_2^2 \frac{\tau(\eta^5)^2+6\tau(\eta^5)+1}{12}-\frac{\mathfrak{r}(\eta)^2+6\mathfrak{r}(\eta)+1}{12},$$
if $(X_2,Y_2)$ are the coordinates on $E_5(\xi^5)$ and $(X_3,Y_3)$ are the coordinates on $E_5(\eta^{5\sigma^{-1}})$.  Therefore, the points with $X_2=0$ map to points with
$$X_3 =\frac{(-5 + \sqrt{5})(\eta \sqrt{5} + 2\eta^2 - \sqrt{5} - 3\eta + 3)(\eta \sqrt{5} - 2\eta^2 - \sqrt{5} + 3\eta - 3)}{20(-2\eta^2 + \sqrt{5} - 2\eta - 3)}.$$
\medskip

Finally, we choose $u=\frac{1}{\varepsilon \xi} \in \textsf{K}_5(\sqrt{5})$, so that
$$u^5 = \frac{1}{\varepsilon^5 \xi^5}=-\bar \varepsilon^5 \frac{\varepsilon^5 \eta^5 +1}{-\eta^5+\varepsilon^5} = - \frac{\eta^5- \bar \varepsilon^5}{\eta^5-\varepsilon^5},$$
as required above for the formula $X(P)$.  Then we compute that
$$u^{\sigma^{-1}} = \frac{1}{\varepsilon \xi^{\sigma^{-1}}}= \frac{1}{\varepsilon T(\eta)},$$
which implies that $\eta=T\left(\varepsilon^{-1} u^{-\sigma^{-1}}\right)$.  Substituting this expression for $\eta$ in $X_3$ gives
$$X_3=\frac{-\varepsilon^4}{2}\frac{(- 2u_1^2+(1+\sqrt{5})u_1 - 3\sqrt{5} - 7)(2u_1^2 + (2\sqrt{5} + 4)u_1 + 3\sqrt{5} + 7)}{(- 2u_1^2+(\sqrt{5}+1)u_1 - 2)(u_1 + 1)^2},$$
with $u_1=u^{\sigma^{-1}}$.  Comparing with the above formula for $X(P)$ shows that $X_3 = X(P)^{\sigma^{-1}}$ and therefore the points $\pm (0,0)$ on $E_5(\xi^5)$ map to $\pm P^{\sigma^{-1}}$ on $E_5(\eta^{5\sigma^{-1}})$. \medskip

This discussion shows that the isogeny $\phi_1 =\iota_2 \circ \iota_1 \circ \psi$ from $E_5(\eta^5)$ to $E_5(\eta^5)^{\sigma^{-1}}$ satisfies
$$\phi_1(P) = \pm P^{\sigma^{-1}}.$$
Applying $\sigma^{-i+1}$ to this gives $\phi_i(P^{\sigma^{-i+1}}) = \pm P^{\sigma^{-i}}$, and therefore
$$\iota(P) = \phi_n \circ \phi_{n-1} \circ \cdots \circ \phi_1(P) = \pm P^{\sigma^{-n}} = \pm P.$$
Since $P$ is a point of order $5$ on $E_5(\eta^5)$, and $P$ does not lie in $\textrm{ker}(\iota)$, we see that $\iota$ is indeed a cyclic isogeny. \medskip

From this and the fact that $\textrm{deg}(\iota)=5^n$ we conclude that the $j$-invariant $j_\eta=j(E_5(\eta^5))$ satisfies the modular equation
$$\Phi_{5^n}(j_\eta,j_\eta)=0.$$
On the other hand, from \cite[p. 263]{co},
$$\Phi_{5^n}(X,X) = c_n \prod_{-d}{H_{-d}(X)^{r(d,5^n)}},$$
where the product is over the discriminants of orders $\textsf{R}_{-d}$ of imaginary quadratic fields and 
$$r(d,5^n)=|\{\alpha \in \textsf{R}_{-d}: \ \alpha \ \textrm{primitive}, \ N(\alpha) = 5^n\}/\textsf{R}_{-d}^{\times}|.$$
Thus, $r(d,5^n)$ is nonzero only when $4^k \cdot 5^n =x^2+dy^2$, ($k=0,1$), has a primitive solution.  Now the polynomial $\textsf{P}_n(x) \in \mathbb{Z}[x]$ splits completely in $\textsf{K}_5(\sqrt{5})$, and its ``unramified'' roots all lie in $\textsf{K}_5$.  Furthermore the ``ramified'' roots all have the form $\xi=T(\eta^\sigma)$ for some unramified root $\eta$, and the corresponding $j$-invariants have the form
$$j_\xi=\frac{(1-12\xi^5+14\xi^{10}+12\xi^{15}+\xi^{20})^3}{\xi^{25}(1-11\xi^5-\xi^{10})},$$
which equals
$$j_\xi  = \frac{(1+228\eta^5+494\eta^{10}-228\eta^{15}+\eta^{20})^3}{\eta^5(1-11\eta^5-\eta^{10})^5}.$$
It follows that all the $j$-invariants $j_\eta, j_\xi$ lie in $\textsf{K}_5$.  Hence, the value $d$ for which $H_{-d}(j_\eta)=0$ is not divisible by $5$.  Thus, $(5,xyd)=1$, and therefore $\left(\frac{-d}{5}\right)=+1$. \medskip

From $H_{-d}(j_\eta) = H_{-d}((j_\eta)^{\sigma^{-1}}) = H_{-d}(j_\xi) = 0$ we see that the periodic point $\eta$ is a root of both polynomials $F_d(x^5), G_d(x^5)$, where
\begin{align*}
F_d(x) = & \ x^{5h(-d)}(1-11x-x^2)^{h(-d)}H_{-d}\left(\frac{(x^4+12x^3+14x^2-12x+1)^3}{x^5(1-11x-x^2)}\right),\\
G_d(x) = & \ x^{h(-d)}(1-11x-x^2)^{5h(-d)}H_{-d}\left(\frac{(x^4-228x^3+494x^2+228x+1)^3}{x(1-11x-x^2)^5}\right).
\end{align*}
Now the roots of the polynomial $G_d(x^5)$ are invariant under the action of the icosahedral group $G_{60}=\langle S, T \rangle$, where $T$ is as before and $S(z) = \zeta z$, with $\zeta = e^{2 \pi i/5}$.  (See \cite{mor}, \cite{mor9}.)  Since $H_{-d}(X)$ is irreducible over the field $L=\mathbb{Q}(\zeta)$, containing the coefficients of all the maps in $G_{60}$, the polynomial $G_d(x^5)$ factors over $L$ into a product of irreducible polynomials of the same degree.  (See the similar argument in \cite[p. 864]{mor4}.)  By the results of \cite[pp. 1193, 1202]{mor8}, one of these irreducible factors is $p_d(x)$, whose degree is $4h(-d)$, and $p_d(x)$ is invariant under the action of the subgroup
$$H = \langle U, T \rangle, \ \ U(z)=\frac{-1}{z},$$
a Klein group of order $4$.  The normalizer of $H$ in $G_{60}$ is $N=\langle A, H \rangle \cong A_4$, where $A=STS^{-2}$ is the map
$$A(z)=\zeta^3 \frac{(1+\zeta)z+1}{z-1-\zeta^4}$$
of order $3$, and $ATA^{-1}=U, AUA^{-1}=T_2=TU$.  The distinct left cosets of $H$ in $G_{60}$ are represented by the elements
$$M_{ij}=S^j A^i, \ \ 0 \le i \le 2, \ 0 \le j \le 4.$$
(See \cite[Prop. 3.3]{mor9}.)  We would like to show that $\eta$ is a root of the factor $p_d(x)$. \medskip

Since all the roots of $G_d(x^5)$ have the form $M_{ij}(\alpha)$, for some root $\alpha$ of $p_d(x)$ (\cite[p. 1203]{mor8}), the factors of $G_d(x^5)$ over $L$ have the form
$$p_{i,j}(x) = (cx+d)^{4h(-d)}p_d(A^iS^j(x)),$$
where $A^i S^j(x) = \frac{ax+b}{cx+d}$.  The stabilizer of this polynomial in $G_{60}$ is
$$(A^iS^j)^{-1}H A^iS^j = S^{-j}H S^j,$$
which contains the map $S^{-j}US^j(x)=\frac{-\zeta^{-2j}}{x}$.  If $p_{i,j}(\eta)=0$, where $j \neq 0$, then both $\eta$ and $\frac{-\zeta^{-2j}}{\eta}$ are roots of $p_{i,j}(x)$, which would imply that $\zeta^{-2j}$ is contained in the splitting field of $\textsf{P}_n(x)$ over $\mathbb{Q}$, and is therefore contained in $\textsf{K}_5(\sqrt{5})$, which is not the case.  Hence, $\eta$ can only be a root of $p_{i,0}(x) = (c_ix+d_i)^{4h(-d)}p_d(A^i(x))$, for some $i$.  But then the elements in $HA^i(\eta)$ are roots of $p_d(x)$.  Assume $i=1$.  Since $A(\eta)$ is a root of $p_d(x)$, so is $A^{\rho^j}(\eta)$, where $\rho$ is the automorphism of $\textsf{K}_5(\zeta)/\textsf{K}_5$ for which $\zeta^\rho=\zeta^2$.  But $A^\rho = A^{-1}U$, so that $A^{\rho^2}=A^{-\rho}U=UAU$ and $A^{\rho^3} = UA^\rho U =UA^{-1}$.  Thus, $A^{\rho^3}(\eta)$ being a root of $p_d(x)$ and $U \in H$ imply that $A^{-1}(\eta)$ is also a root of $p_d(x)$.  But then $\eta$ is a common root of $p_{1,0}(x) = (c_1x+d_1)^{4h(-d)}p_d(A(x))$ and $p_{2,0}(x)=(c_2x+d_2)^{4h(-d)}p_d(A^{-1}(x))$, which is impossible, since these are two of the irreducible factors of $G_d(x^5)$ over $L$, and the latter polynomial has no multiple roots, for $d \neq 4$.  (See \cite[\S 2.2]{mor9}.)  A similar argument works if $i=2$, since $A^2 = A^{-1}$ and $A = UA^{-\rho}$.  For $d=4$, we have
\begin{align*}
G_4(x^5)= & \ (x^{20}-228x^{15}+494x^{10}+228x^5+1)^3-1728 x^5(1-11x^5-x^{10})^5\\
 = & \ (x^2 + 1)^2(x^4 + 2x^3 - 6x^2 - 2x + 1)^2(x^8 - x^6 + x^4 - x^2 + 1)^2\\
& \times \ (x^8 + 4x^7 + 17x^6 + 22x^5 + 5x^4 - 22x^3 + 17x^2 - 4x + 1)^2\\
& \times \ (x^8 - 6x^7 + 17x^6 - 18x^5 + 25x^4 + 18x^3 + 17x^2 + 6x + 1)^2,
 \end{align*}
and the only periodic point $\eta \in \textsf{D}_5$ which is a root of $G_4(x^5)$ is the fixed point
$$\eta= i = 3+3 \cdot 5+2 \cdot 5^2+3 \cdot 5^3+ 5^4 + \cdots \ \ \in \mathbb{Q}_5.$$
Thus, $d=4$ does not occur when $n \ge 2$.  (Except for the primitive $20$-th roots of unity, which do not lie in $\textsf{K}_5(\sqrt{5})$, the other roots of $G_4(x^5) = 0$ satisfy $x \equiv 2$ mod $5$, and so do not lie in $\textsf{D}_5$.) \medskip

Hence, the only possibility is that $p_d(\eta)=0$.  This shows that all periodic points of $T_5(x)$ in $\textsf{D}_5$ are roots of some $p_d(x)$ for which $(-d/5)=+1$.  Since $T_5(\eta) = \eta^{\tau_5}$ for such a root by (\ref{eq:3.6}), it is clear that $\tau_5$ has order $n$ in the corresponding Galois group $\textrm{Gal}(F_1/\mathbb{Q})$, as well.  All the roots of $\textsf{P}_n(x)$ which do not lie in $\textsf{D}_5$ have the form $T(\eta)$, for $\eta \in \textsf{D}_5$, by the discussion in Section 2, and are also roots of $p_d(x)$ for one of these integers $d$, since $T(x)$ stabilizes the roots of $p_d(x)$.  \medskip 

Thus, if $n \ge 2$, the only irreducible factors of $\textsf{P}_n(x)$ over $\mathbb{Q}$ are the polynomials $p_d(x)$ for which $(-d/5) = +1$ and $\tau_5 \in \textrm{Gal}(F_1/\mathbb{Q})$ has order $n$.  This proves (\ref{eq:4.1}).  $\square$  \medskip

For use in the following corollary, note that the substitution $(X,Y) \rightarrow \left(\frac{-1}{X},\frac{-1}{Y}\right)$ represents an automorphism of the curve $g(X,Y)=0$, since
\begin{equation}
X^5 Y^5 g\left(\frac{-1}{X},\frac{-1}{Y}\right)=g(X,Y).
\label{eq:4.2}
\end{equation}
As in \cite{mor8}, put
\begin{equation}
g_1(X,Y)=Y^5 g\left(X,\frac{-1}{Y}\right).
\label{eq:4.3}
\end{equation}

In the following corollary, we prove the claim stated in the last paragraph of \cite[p. 1212]{mor8}.  In that paragraph, the polynomial $x^2+x-1$ should have also been listed along with $x, x^2+1$ and $p_d(x)$ as factors of the resultants $R_n(x)$.  As we will see below, however, $x^2+x-1$ never divides $\tilde R_n(x)$. \bigskip

\noindent {\bf Corollary 4.3.} {\it If $\alpha \neq 0$ is a root of the $(n-1)$-fold iterated resultant
$$\tilde R_n(x) = Res_{x_{n-1}}(...(Res_{x_2}(Res_{x_1}(g(x,x_1),g(x_1,x_2)),g(x_2,x_3)),...,g_1(x_{n-1},x)),$$
for $n \ge 2$, then $\alpha$ is either $\pm i$ or a root of some polynomial $p_d(x)$, where $p_d(x) \mid R_{2n}(x)$.} \smallskip

\noindent {\it Proof.} A root $\alpha \neq 0$ of $\tilde R_n(x)$ satisfies the simultaneous equations
$$g(\alpha, \alpha_1) = g(\alpha_1, \alpha_2) = \cdots = g(\alpha_{n-2}, \alpha_{n-1}) = g_1(\alpha_{n-1}, \alpha) = 0,$$
for some elements $\alpha_i$ in $\overline{\mathbb{Q}}$, the algebraic closure of $\mathbb{Q}$.  Note that $\alpha_i \neq 0$, for $1 \le i \le n-1$, because $g(X,0) = X^5$, so that $\alpha_i = 0$ implies $\alpha_{i-1}=0$.  But the definition of $g_1(X,Y)$ and the final equation in the above chain give that $g\left(\alpha_{n-1}, \frac{-1}{\alpha}\right) = 0$.  Now the identity (\ref{eq:4.2}) implies, using the above simultaneous equations, that
$$g\left(\frac{-1}{\alpha}, \frac{-1}{\alpha_1}\right) = g\left(\frac{-1}{\alpha_1},\frac{-1}{\alpha_2}\right) = \cdots = g\left(\frac{-1}{\alpha_{n-1}}, \alpha\right) = 0.$$
Tacking this chain of equations onto the first chain following the equation $g\left(\alpha_{n-1}, \frac{-1}{\alpha}\right) = 0$ shows that $\alpha$ is a root of $R_{2n}(x)=0$.  Setting $p_4(x) = x^2+1$ (see below), we only have to verify that $\alpha$ is not a root of $x^2+x-1$ to conclude that $\alpha$ is a root of some polynomial $p_d(x)$, because 
$$\textsf{P}_1(x) = x(x^2 + 1)(x^2 + x - 1)(x^4 + x^3 + 3x^2 - x + 1) = x(x^2+x-1) p_4(x) p_{19}(x).$$
For in that case $\alpha$ is either a root of $p_4(x) p_{19}(x)$ or a root of some $\textsf{P}_m(x)$, for $m > 1$.  But if $\alpha = \frac{-1 \pm \sqrt{5}}{2}$, then $\alpha$ is a fixed point, $g(\alpha,y) = 0 \Rightarrow y = \alpha$, but
$$g_1(\alpha,\alpha) = \alpha^5 g(\alpha, \overline{\alpha}) = \frac{625 - 275 \sqrt{5}}{2} \neq 0.$$
Thus, $\alpha$ cannot be a root of $\tilde R_n(x)$ for any $n \ge 1$.   $\square$  \medskip

\noindent {\bf Remark.} This justifies the claims made in Section 5 of Part II about the resultant $\tilde R_n(x)$.  In particular, all its irreducible factors are $x^2+1$ and polynomials of the form $p_d(x)$.  This shows also that the polynomial in Example 2 of that section (pp. 1210-1211) is indeed $p_{491}(x)$.  The computation of the degree $\tilde R_3(x)$ was in error, however, at the beginning of that example.  In fact the degree is $250$, and there are five factors of degree $12$, not three, as was claimed before: these factors are the polynomials $p_d(x)$ for $d= 31, 44, 124, 211, 331$.  \bigskip

Note that the root $-i =r\left(\frac{-7+i}{5}\right)$, so $p_4(x)$ is the minimal polynomial of a value $r(w_4/5)$, with $w_4=-7+i \in \mathbb{Q}(\sqrt{-4})$ and $\wp_5^2 =(-2+i)^2 \mid w_4$.  This justifies the notation $p_4(x)$. See \cite[p. 139]{du}.  \medskip

The following theorem is immediate from Theorem 4.2 and the computations of Section 2. \bigskip

\noindent {\bf Theorem 4.4.}  {\it The set of periodic points in $\overline{\mathbb{Q}}$ (or $\overline{\mathbb{Q}}_5$ or $\mathbb{C}$) of the multi-valued algebraic function $\mathfrak{g}(z)$ defined by the equation $g(z,\mathfrak{g}(z))=0$ consists of $0, \frac{-1\pm \sqrt{5}}{2}$, and the roots of the polynomials $p_d(x)$, for negative discriminants $-d$ satisfying $\left(\frac{-d}{5}\right)=+1$.  Over $\overline{\mathbb{Q}}$ or $\mathbb{C}$ the latter values coincide with the values $\eta = r(w_d/5)$ and their conjugates over $\mathbb{Q}$, where $r(\tau)$ is the Rogers-Ramanujan continued fraction and the argument $w_d \in K=\mathbb{Q}(\sqrt{-d})$ satisfies}
$$w_d=\frac{v+\sqrt{-d}}{2} \in R_K, \ \ \wp_5^2 \mid w_d, \ \ \textrm{and} \ (N(w_d),f) = 1.$$

The fixed points $0, \frac{-1\pm \sqrt{5}}{2}$ come from the factors $x, x^2+x-1$ of the polynomial $\textsf{P}_1(x)$.  \medskip

Equating degrees in the formula (\ref{eq:4.1}) yields
$$\textrm{deg}(\textsf{P}_n(x))=\sum_{-d \in \mathfrak{D}_{n,5}}{4h(-d)}, \ \ n >1.$$
From (\ref{eq:2.7}) we get the following class number formula. \bigskip

\noindent {\bf Theorem 4.5.} {\it For $n > 1$ we have
$$\sum_{-d \in \mathfrak{D}_{n,5}}{h(-d)} = \frac{1}{2} \sum_{k \mid n}{\mu(n/k)5^k},$$
where $ \mathfrak{D}_{n,5}$ has the meaning given in Theorem 1.3.} \bigskip

Note that the corresponding formula for $n=1$ reads
$$\sum_{-d \in \mathfrak{D}_{1,5}}{h(-d)} = h(-4) + h(-19) = 2 = \frac{1}{2}(5-1).$$

\section{Ramanujan's modular equations for $r(\tau)$}

In this section we take a slight detour to show how the polynomials $p_{4d}(x), p_{9d}(x)$ and $p_{49d}(x)$ can be computed, if the polynomial $p_d(x)$ is known.  \medskip

From Berndt's book \cite[p. 17]{ber} we take the following identity relating $u=r(\tau)$ and $v=r(3\tau)$:
\begin{equation}
(v-u^3)(1+uv^3)=3u^2v^2.
\label{eq:5.1}
\end{equation}
Let
$$P_3(u,v)=(v-u^3)(1+uv^3)-3u^2v^2.$$
This polynomial satisfies the identity
$$v^4 P_3\left(u,\frac{-1}{v}\right) = P_3(v,u).$$
The following theorem gives a simple method of calculating $p_{9d}(x)$ from $p_d(x)$. \bigskip

\noindent {\bf Theorem 5.1.} {\it For any negative discriminant $-d \equiv \pm 1$ (mod $5$), the polynomial $p_{9d}(x)$ divides the resultant}
$$\textrm{Res}_y(P_3(y,x),p_d(y)).$$

\noindent {\it Proof.}  Let $-d=d_Kf^2$, where $d_K$ is the discriminant of $K=\mathbb{Q}(\sqrt{-d})$.  One of the roots of $p_{9d}(x)$ is $\eta' = r(w_{9d}/5)$, where $w_{9d} =\frac{v+\sqrt{-9d}}{2} \in \textsf{R}_{-9d}$, $\wp_5^2 \mid w_{9d}$ and $N(w_{9d})=\frac{v^2+9d}{4}$ is prime to $3f$.  Let $f=3^s f'$, with $(f',3)=1$.  For some integer $k$, $w_{9d}+25f'k = \frac{v+50f'k+\sqrt{-9d}}{2}$ satisfies $v+50f'k \equiv v-4f'k \equiv 3$ mod $9$.  Furthermore,
$$\eta'=r\left(\frac{w_{9d}+25f'k}{5}\right) = r\left(\frac{w_{9d}}{5}+5f'k\right) = r\left(\frac{w_{9d}}{5}\right).$$
Thus, we may assume $3 \ || \ v$, and then $9 \mid N(w_{9d})$.  In that case $w_d=\frac{w_{9d}}{3} \in \textsf{R}_{-d}$, where $(N(w_d),f)=1$, even when $3 \mid f$.  Furthermore, $\wp_5^2 \mid w_d$.  Hence, $\eta = r(w_d/5)$ is a root of $p_d(x)$.  From (\ref{eq:5.1}) we have
$$P_3(\eta,\eta')=P_3(r(w_d/5),r(w_{9d}/5))=P_3(r(w_d/5),r(3w_d/5))=0.$$
Hence, $\eta'$ is a root of the resultant, which therefore has its minimal polynomial $p_{9d}(x)$ as a factor. $\square$ \bigskip

\noindent {\bf Example 1.} We compute
$$\textrm{Res}_y(P_3(y,x),p_4(y)) = \textrm{Res}_y(P_3(y,x),y^2+1) = x^8 + x^6 - 6x^5 + 9x^4 + 6x^3 + x^2 + 1.$$
Since the latter polynomial is irreducible, the theorem shows that it equals $p_{36}(x)$:
$$p_{36}(x)=x^8 + x^6 - 6x^5 + 9x^4 + 6x^3 + x^2 + 1.$$
This verifies once again the entry for $d=36$ in Table 1 of \cite{mor8}, which we used in Example 1 of that paper (p. 1208).  In the same way, we compute
\begin{align*}
\textrm{Res}_y(P_3(y,x), & \ p_{36}(y)) = (x^2 + 1)^4(x^{24} - 18x^{23} + 81x^{22} - 60x^{21} + 594x^{20} \\
& + 1074x^{19} + 118x^{18} - 1002x^{17} - 261x^{16} + 6882x^{15} + 12078x^{14} \\ 
& + 1014x^{13} - 18585x^{12} - 1014x^{11} + 12078x^{10} - 6882x^9 - 261x^8  \\
& + 1002x^7 + 118x^6 - 1074x^5 + 594x^4 + 60x^3 + 81x^2 + 18x + 1)\\
& = p_4(x)^4 p_{324}(x).
\end{align*}

There is also the identity from \cite[p. 12]{ber} relating $u=r(\tau)$ and $v=r(2\tau)$:
\begin{equation}
(v-u^2)=(v+u^2) \cdot uv^2.
\label{eq:5.2}
\end{equation}
Setting
$$P_2(u,v)=(v+u^2) \cdot uv^2-(v-u^2),$$
we have the following identity, analogous to the identity for $P_3(u,v)$. 
$$v^3 P_2\left(u,\frac{-1}{v}\right) = P_2(v,u).$$
An argument similar to the proof of Theorem 5.1 yields \bigskip

\noindent {\bf Theorem 5.2.}  {\it For any negative discriminant $-d \equiv \pm 1$ (mod $5$), the polynomial $p_{4d}(x)$ divides the resultant}
$$\textrm{Res}_y(P_2(y,x),p_d(y)).$$
\noindent {\it Proof.}  Again, let $-d=d_Kf^2$, where $d_K$ is the discriminant of $K=\mathbb{Q}(\sqrt{-d})$.  One of the roots of $p_{4d}(x)$ is $\eta' = r(w_{4d}/5)$, where $w_{4d} =\frac{v+\sqrt{-4d}}{2} \in \textsf{R}_{-4d}$, $\wp_5^2 \mid w_{4d}$ and $N(w_{4d})=\frac{v^2+4d}{4}$ is prime to $2f$.  Thus, $v \equiv 2d +2$ (mod $4$).  If $f$ is odd, we set
$$w'=w_{4d}+25f = \left(\frac{v}{2}+25f\right)+\sqrt{-d}=v'+\sqrt{-d}.$$
Then,
$$r\left(\frac{w'}{5}\right) = r\left(\frac{w_{4d}}{5}+5f\right) = r\left(\frac{w_{4d}}{5}\right) =\eta'.$$
Moreover, $v' \equiv \frac{v}{2}+1 \equiv d$ (mod $2$).  Now let $w_d=\frac{w'}{2}=\frac{v'+\sqrt{-d}}{2} \in \textsf{R}_{-d}$, where $(N(w_d),f)=1$.  Then $\wp_5^2 \mid w_d$ and $\eta = r(w_d/5)$ is a root of $p_d(x)$.  From (\ref{eq:5.2}) we have
$$P_2(\eta,\eta')=P_2(r(w_d/5),r(w_{4d}/5))=P_2(r(w_d/5),r(2w_d/5))=0.$$
Hence, $\eta'$ is a root of the resultant, which therefore has its minimal polynomial $p_{4d}(x)$ as a factor. \medskip

On the other hand, if $f$ is even, let $f=2^s f'$, with $f'$ odd.  Then $d$ is even, so $v/2$ is odd.  In this case we choose $k$ so that
$$v' = \frac{v}{2}+25f'k \equiv \begin{cases} 0 \ (\textrm{mod} \ 4), & \textrm{if} \ 4 \ || \ d;\\ 2  \ (\textrm{mod} \ 4), & \textrm{if} \ 8 \mid d. \end{cases}$$
With this choice of $k$ we have $v' \equiv d$ (mod $2$), so letting $w' = v'+\sqrt{-d}=w_{4d}+25f'k$ and $w_d = \frac{w'}{2}$, we have
$w_d  \in \textsf{R}_{-d}$ and
$$N(w_d) = \frac{v'^2+d}{4} \equiv \begin{cases} \frac{d}{4} \equiv 1 \ (\textrm{mod} \ $2$), & \textrm{if} \ 4 \ || \ d;\\ \frac{v'^2}{4} \equiv 1 \ (\textrm{mod} \ $2$), & \textrm{if} \ 8 \mid d. \end{cases}$$
In either case, we get that $(N(w_d),f) = 1$.  We have $r(w'/5)=r(w_{4d}/5)$, as before, and letting $\eta=r(w_d/5)$ be a root of $p_d(x)$, we obtain $P_2(\eta,\eta') = 0$ as above, and the assertion of the theorem follows.  $\square$ \bigskip

\noindent {\bf Example 2.} We have
\begin{align*}
\textrm{Res}_y(P_2(y,x), & \ p_{36}(y)) = (x^8 + x^6 - 6x^5 + 9x^4 + 6x^3 + x^2 + 1)\\
& \times (x^{16} - 2x^{15} + 18x^{14} + 24x^{13} + 83x^{12} + 78x^{11} + 74x^{10} + 40x^9\\
& \ \ + 9x^8 - 40x^7 + 74x^6 - 78x^5 + 83x^4 - 24x^3 + 18x^2 + 2x + 1)\\
& = p_{36}(x) p_{144}(x)
\end{align*}
and
\begin{align*}
\textrm{Res}_y&(P_2(y,x), \ p_{144}(y)) = (x^8 + x^6 - 6x^5 + 9x^4 + 6x^3 + x^2 + 1)^2\\
& \times \ (x^{32} - 32x^{31} + 586x^{30} - 2856x^{29} + 5818x^{28} - 160x^{27} - 23408x^{26}\\
& \ + 41964x^{25} - 6573x^{24} - 63520x^{23} + 64426x^{22} + 12736x^{21} - 38746x^{20}\\
& \ - 11464x^{19} + 55416x^{18} - 38148x^{17} - 5743x^{16} + 38148x^{15} + 55416x^{14}\\
& \ + 11464x^{13} - 38746x^{12} - 12736x^{11} + 64426x^{10} + 63520x^9 - 6573x^8\\
& \ - 41964x^7 - 23408x^6 + 160x^5 + 5818x^4 + 2856x^3 + 586x^2 + 32x + 1)\\
& = p_{36}(x)^2 p_{576}(x).
\end{align*}

We can use Theorems 5.1 and 5.2 to construct polynomials $p_d(x)$ for which the Conjecture (1) in \cite[p. 1199]{mor8} does not hold.  For example, starting with
$$p_{51}(x) = x^8 + x^7 + x^6 - 7x^5 + 12x^4 + 7x^3 + x^2 - x + 1,$$
applying Theorem 5.2 once gives that
\begin{align*}
p_{204}(x) = & \ x^{24} - x^{23} + 38x^{22} + 36x^{21} + 166x^{20} + 33x^{19} + 57x^{18} + 22x^{17}\\
& + 573x^{16} + 1603x^{15} + 2465x^{14} + 1225x^{13} + 1768x^{12} - 1225x^{11}\\
& + 2465x^{10} - 1603x^9 + 573x^8 - 22x^7 + 57x^6 - 33x^5 + 166x^4 - 36x^3\\
& + 38x^2 + x + 1,
\end{align*}
whose discriminant is exactly divisible by $17^{12}$, in accordance with Conjecture (1).  Applying Theorem 5.2 to this polynomial yields
the polynomial $p_{816}(x)$, of degree $48$, whose discriminant is exactly divisible by $17^{40}$:
$$\textrm{disc}(p_{816}(x)) = 2^{160}3^{120}5^{276}7^{40}17^{40}31^{24}47^8 79^8 179^4 191^{12}241^8 491^8 541^8 691^8;$$
whereas Conjecture (1) predicts that $17^{24}$ should be the power of $17$ dividing $\textrm{disc}(p_{816}(x))$. \medskip

Note that the period of the roots of $p_{51}(x)$ is $4$, whereas the period of the roots of $p_{204}(x)$ and $p_{816}(x)$ is $12$.
\medskip

We modify the statement of Conjecture (1) in \cite[p. 1199]{mor8} as follows. \bigskip

\noindent {\bf Conjecture (1$'$).} {\it If $q >5$ is a prime which divides the field discriminant $d_K$ of $K=\mathbb{Q}(\sqrt{-d})$, then $q^{2h(-d_K)}$ exactly divides $\textrm{disc}(p_{d_K}(x))$.}  \bigskip

Now define the polynomial $P_7(u,v)$ by
\begin{align*}
P_7(u,v) = & \ u^8 v^7+(-7v^5+1)u^7+7u^6v^3+7(-v^6+v)u^5+35u^4 v^4\\
& + 7(v^7+v^2)u^3-7u^2v^5-(v^8+7v^3)u-v.
\end{align*}
Note that $P_7(u,v)$ satisfies the polynomial identity
$$v^8P_7\left(u,\frac{-1}{v}\right) = P_7(v,u).$$
We make the \bigskip

\noindent {\bf Conjecture.} {\it The Rogers-Ramanujan function $r(\tau)$ satisfies
$$P_7(r(\tau), r(7\tau)) = 0.$$}
A calculation on Maple shows that the first $250$ $q$-coefficients of $P_7(r(\tau), r(7\tau))$ are zero. \bigskip

\noindent {\bf Theorem 5.3.} {\it This conjecture implies that for any negative discriminant $-d \equiv \pm 1$ (mod $5$), the polynomial $p_{49d}(x)$ divides the resultant}
$$\textrm{Res}_y(P_7(y,x),p_d(y)).$$

The proof is the same, mutatis mutandis, as the proof of Theorem 5.1, on replacing the prime $3$ by $7$. \medskip

\noindent {\bf Example 3.} We compute that
\begin{align*}
\textrm{Res}_y(P_7(y,x), p_4(y)) = & \ p_{196}(x)\\
= & \ x^{16} + 14 x^{15} + 64 x^{14} + 84 x^{13} - 35 x^{12} - 14 x^{11} + 196 x^{10}\\
& \ + 672 x^9 + 1029 x^8 - 672 x^7 + 196 x^6 + 14 x^5 - 35 x^4\\
& \ - 84 x^3 + 64 x^2 - 14 x + 1.
\end{align*}
As a check, note that $h(-4 \cdot 7^2) = 4$ and the discriminant of $p_{196}(x)$ is
$$\textrm{disc}(p_{196}(x)) = 2^{32} \cdot 3^{12} \cdot 5^{28} \cdot 7^{14} \cdot 19^4 \cdot 71^8,$$
all of whose prime factors are less than $d=196 = 4 \cdot 7^2$.

\section{Periodic points for $h(t,u)$.}

\subsection{Reduction to periodic points of $g(x,y)$.}

From \cite{mor8} the equation connecting $t=X-\frac{1}{X}$ and $u=Y-\frac{1}{Y}$ in the function field of the curve $g(X, Y) = 0$ is
\begin{align*}
h(t,u)=u^5&-(6+5t+5t^3+t^5)u^4+(21+5t+5t^3+t^5)u^3\\
&-(56+30t+30t^3+6t^5)u^2+(71+30t+30t^3+6t^5)u\\
&-120-55t-55t^3-11t^5.
\end{align*}
On this curve $\upsilon = \eta - \frac{1}{\eta} \in \Omega_f$, with $\eta = r(w_d/5)$, satisfies
$$h(\upsilon, \upsilon^{\tau_5}) = 0, \ \ \ \tau_5 = \left(\frac{\Omega_f/\mathbb{Q}(\sqrt{-d})}{\wp_5}\right).$$
This yielded the following theorem. \bigskip

\noindent {\bf Theorem 6.1.}  {\it If $\Omega_f$ is the ring class field of conductor $f$ (relatively prime to $5$) over the field $K=\mathbb{Q}(\sqrt{-d})$, where $-d=d_Kf^2$ and $\left(\frac{-d}{5}\right)=+1$, then $\Omega_f=K(\upsilon)$, where $\upsilon=\eta-\frac{1}{\eta}$ is a periodic point of the algebraic function $\mathfrak{f}(z)$ defined by $h(z,\mathfrak{f}(z))=0$.} \bigskip

Note the identity
\begin{equation}
X^5 Y^5 h\left(X-\frac{1}{X},Y-\frac{1}{Y}\right) = -g(X,Y)g_1(X,Y),
\label{eq:6.1}
\end{equation}
where $g(X,Y)$ is given by (\ref{eq:2.1}) and $g_1(X,Y)$ is defined in (\ref{eq:4.3}).  Also, recall that
\begin{equation}
X^5 Y^5 g\left(\frac{-1}{X},\frac{-1}{Y}\right) = g(X,Y), \ \ X^5 Y^5 g_1\left(\frac{-1}{X},\frac{-1}{Y}\right) = g_1(X,Y),
\label{eq:6.2}
\end{equation}
where the second identity is an easy consequence of the first.  Using these facts we can prove the following. \bigskip

\noindent {\bf Theorem 6.2.} {\it If $\upsilon \neq -1$ is any periodic point of the algebraic function $\mathfrak{f}(z)$ in Theorem 6.1, then
$$\upsilon  = \eta - \frac{1}{\eta},$$
for some periodic point $\eta$ of $\mathfrak{g}(z)$, and $\upsilon$ generates a ring class field $\Omega_f$ over some field $K = \mathbb{Q}(\sqrt{-d})$, where $-d=d_Kf^2$ and $\left(\frac{-d}{5}\right) = +1$.} \medskip

\noindent {\it Proof.} Assume that there exist elements $\upsilon_i$ for which
\begin{equation}
h(\upsilon,\upsilon_1) = h(\upsilon_1,\upsilon_2) = \cdots = h(\upsilon_{n-1},\upsilon) = 0.
\label{eq:6.3}
\end{equation}
Since the substitution $x=y-\frac{1}{y}$ transforms the polynomial
$$h(x,x) = -(x + 1)(x^2 + 4)(x^2 - x + 3)(x^2 - 2x + 2)(x^2 + x + 5),$$
(after multiplying by $y^9$) into the product
\begin{align*}
-(y^2 + y - 1)&(y^2 + 1)^2(y^4 - y^3 + y^2 + y + 1)(y^4 - 2y^3 + 2y + 1)\\
& \ \ \ \ \times (y^4 + y^3 + 3y^2 - y + 1)\\
& =-(y^2 + y - 1) p_4(y)^2 p_{11}(y) p_{16}(y) p_{19}(y),
\end{align*}
we may assume $n \ge 2$.  Set $g_0(X,Y) = g(X,Y)$ and write $\upsilon=\eta-\frac{1}{\eta}$ and $\upsilon_i = \eta_i-\frac{1}{\eta_i}$.  By (\ref{eq:6.1}), equation (\ref{eq:6.3}) is equivalent to a set of simultaneous equations
\begin{equation}
g_{i_1}(\eta,\eta_1)=g_{i_2}(\eta_1,\eta_2) = \cdots = g_{i_n}(\eta_{n-1},\eta) = 0,
\label{eq:6.4}
\end{equation}
where each $i_k = 0$ or $1$.  Using the same idea as in the proof of Corollary 4.3, we will transform this set of equations into a set of equations which only involve the polynomial $g=g_0$.  Assume first that $i_1=1$.  Then
$$0=g_1(\eta, \eta_1) = g\left(\eta, \frac{-1}{\eta_1}\right).$$
Now we use (\ref{eq:6.2}) to rewrite the remaining equations, so that we have
$$0 = g\left(\eta, \frac{-1}{\eta_1}\right) = g_{i_2}\left(\frac{-1}{\eta_1}, \frac{-1}{\eta_2}\right) = \cdots = g_{i_n}\left(\frac{-1}{\eta_{n-1}}, \frac{-1}{\eta}\right),$$
with the same subscripts $i_r$, for $r \ge 2$, as before.  Now assume we have transformed the first $k-1$ equations so that only the polynomial $g(X,Y)$ appears.  Then, on renaming the elements $\pm \eta_i^{\pm 1}$ as $\eta_i$, we have the simultaneous equations
$$0 = g\left(\eta, \eta_1\right) = \cdots =g(\eta_{k-2},\eta_{k-1}) = g_{i_k}(\eta_{k-1}, \eta_k) = \cdots = g_{i_n}(\eta_{n-1}, \pm \eta^{\pm 1}).$$
If $i_k=0$ we replace $k$ by $k+1$ and continue.  If $i_k=1$ we replace $g_{i_k}(\eta_{k-1}, \eta_k)$ by $g(\eta_{k-1},-1/\eta_k)$ and use (\ref{eq:6.2}) to replace $\eta_r$ in the remaining equations by $-1/\eta_r, r \ge k$.  Then, on renaming the $\eta$'s again, we get a chain of equations
$$0 = g\left(\eta, \eta_1\right) = \cdots =g(\eta_{k-1},\eta_k) =  \cdots = g_{i_n}(\eta_{n-1}, \pm \eta^{\pm 1}).$$
Thus, by induction, we see that (\ref{eq:6.4}) is equivalent to a chain of equations
$$0 = g\left(\eta, \eta_1\right) = \cdots = g(\eta_{n-1}, \pm \eta^{\pm 1})$$
only involving the polynomial $g$.  If the final $\eta$ is simply $\eta$, then $\eta$ is a periodic point of $g$ having period $n$.  On the other hand, if the final $\eta$ appearing in these equations is $-\eta^{-1}$, then we use the same argument as in Corollary 4.3 to show that $\eta$ is a periodic point of period $2n$.  Then we know $\eta$ is not $0$ or a root of $x^2+x-1$, and therefore must be a root of some $p_d(x)$.  By Theorem 6.1, this implies that $K(\upsilon) = \Omega_f$, for $K = \mathbb{Q}(\sqrt{-d})$ and $-d=d_K f^2$.  This proves the theorem.  $\square$ \medskip

Taken together, Theorems 6.1 and 6.2 verify Conjecture 1(b) of Part I for the case $p=5$.  To verify Conjecture 1(a), we define the function
$$\textsf{T}_5(z) = T_5(\eta)-\frac{1}{T_5(\eta)}, \ \ \eta = \frac{z \pm \sqrt{z^2+4}}{2}.$$
We can also write
$$\textsf{T}_5(z) = \phi \circ T_5 \circ \phi^{-1}(z), \ \ \phi(z) = z - \frac{1}{z},$$
where $\phi^{-1}(z) \in  \{\frac{z \pm \sqrt{z^2+4}}{2}\}$ is two-valued.  Since
$$g(z,T_5(z)) = 0 \ \Rightarrow \ g\left(\frac{-1}{z},\frac{-1}{T_5(z)}\right) = 0,$$
it follows from Proposition 3.2 that
$$T_5\left(\frac{-1}{z}\right) = \frac{-1}{T_5(z)}, \ \textrm{for} \ z \in \textsf{D}_5 \cap \{z: |z|_5 = 1\}.$$
Since the two solutions $\eta^{(+)}, \eta^{(-)}$ of $\phi(\eta^{(\pm)})=z$ satisfy $\eta^{(+)} \eta^{(-)} = -1$, the value taken for $\phi^{-1}(z)$ does not affect the value of $\textsf{T}_5(z)$.  In other words, we have the symmetric formula
$$\textsf{T}_5(z) = T_5(\eta^{(+)}) +  T_5(\eta^{(-)}), \ \ \eta^{(\pm)} = \frac{z \pm \sqrt{z^2+4}}{2}.$$
Then from $T_5(\eta^{(+)}) \cdot T_5(\eta^{(-)}) = -1$ and (\ref{eq:3.3}) it follows that $\textsf{T}_5(z) \in \phi(\textsf{D}_5 \cap  \{z: |z|_5 = 1\})$, which implies that
$$\textsf{T}_5^n(z) = T_5^n(\eta^{(+)}) +  T_5^n(\eta^{(-)}), \ n \ge 1, \ \eta^{(\pm)} = \frac{z \pm \sqrt{z^2+4}}{2}.$$
Furthermore, $g(z,T_5(z))=0$ implies that
$$h(z-1/z,\textsf{T}_5(z-1/z)) = -g(z,T_5(z))g_1(z,T_5(z))=0.$$
We deduce the following. \bigskip

\noindent {\bf Theorem 6.3.} {\it For any negative discriminant $-d=d_Kf^2$ with $\left(\frac{-d}{5}\right)=+1$, and for $\eta=r(w_d/5)$, as in Part II, the $h(-d)$ distinct conjugate values
$$\upsilon^\tau = \eta^\tau-\frac{1}{\eta^\tau}, \ \ \tau \in \textrm{Gal}(F_1/K),$$
lying in the ring class field $\Omega_f$ of $K=\mathbb{Q}(\sqrt{-d})$, are periodic points of the $5$-adic algebraic function $\textsf{T}_5(z)$ in the $5$-adic domain
$$\widetilde{\textsf{D}}_5=\phi(\textsf{D}_5 \cap  \{z \in \textsf{K}_5: |z|_5 = 1\}).$$
The period of $\upsilon^\tau$ is equal to the order of the automorphism $\tilde{\tau}_5 =\left(\frac{\Omega_f/K}{\wp_5}\right)$.} \smallskip

\noindent {\it Proof.} This is immediate from
$$\textsf{T}_5(\upsilon^\tau) = \textsf{T}_5\left(\eta^\tau-\frac{1}{\eta^\tau}\right) = T_5(\eta^\tau)-\frac{1}{T_5(\eta^\tau)} = \eta^{\tau \tau_5}-\frac{1}{\eta^{\tau \tau_5}} = \upsilon^{\tau \tau_5},$$
where the third equality above follows from $g(\eta^\tau,\eta^{\tau \tau_5})=0$. The fact that the period is the order of $\tilde \tau$ is a consequence of the fact that $\mathbb{Q}(\upsilon)=\Omega_f$ and that
$$\tilde{\tau}_5 =\tau_5 |_{\Omega_f}, \ \ \tau_5=\left(\frac{F_1/K}{\wp_5}\right).$$
 $\square$  \bigskip

\noindent {\bf Corollary 6.4.} {\it Conjecture 1(a) of \cite{mor5} holds for the prime $p=5$: Every ring class field $\Omega_f$ over $K=\mathbb{Q}(\sqrt{-d})$, with $\left(\frac{-d}{5}\right)=+1$ and $(f,5)=1$, is generated over $\mathbb{Q}$ by a periodic point of the $5$-adic algebraic function $\textsf{T}_5(z)$ which is contained in the domain $\widetilde{\textsf{D}}_5=\phi(\textsf{D}_5 \cap  \{z \in \textsf{K}_5: |z|_5 = 1\}) \subset \textsf{K}_5$.} \bigskip

Note: it is clear that $\textsf{T}_5(\widetilde{\textsf{D}}_5) \subseteq \widetilde{\textsf{D}}_5$, since $T_5(x)$ maps the set $\textsf{D}_5 \cap  \{z \in \textsf{K}_5: |z|_5 = 1\}$ into itself, by Corollary 3.3 and equation (\ref{eq:3.3}). \medskip

The values $\upsilon^\tau$ and their complex conjugates coincide with the roots of the polynomial $t_d(x)$, for which
\begin{equation}
x^{2h(-d)}t_d\left(x-\frac{1}{x}\right)=p_d(x), \ \ d > 4.
\label{eq:6.5}
\end{equation}
Theorem 6.2 shows that every periodic point $\upsilon \neq -1$ of $\mathfrak{f}(z)$ is a root of some polynomial $t_d(x)$.

\subsection{Deuring's class number formula.}

Let
$$S^{(1)}(t,t_1) := h(t,t_1) \equiv 4(t_1+1)^4(t^5-t_1)  \ \ (\textrm{mod} \ 5)$$
and
$$S^{(n)}(t,t_n):=\textrm{Resultant}_{t_{n-1}}(S^{(n-1)}(t,t_{n-1}),h(t_{n-1},t_n)), \ \ n \ge 2.$$
Then it follows by induction that
$$S^{(n)}(t,t_n) \equiv 4(t_n+1)^{5^n-1}(t^{5^n}-t_n)  \ \ (\textrm{mod} \ 5), \ \ n \ge 1.$$
Hence, the polynomial $S_n(t) := S^{(n)}(t,t)$ satisfies the congruence
\begin{equation}
S_n(t) \equiv 4(t+1)^{5^n-1}(t^{5^n}-t)  \ \ (\textrm{mod} \ 5).
\label{eq:6.6}
\end{equation}
It follows that
$$\textrm{deg}(S_n(t))=2 \cdot 5^n-1, \ \ n \ge 1.$$
See the Lemma on pp. 727-728 of Part I, \cite{mor5}. \medskip

Let $L(z)=\frac{-z+4}{z+1}$.  Then
$$L\left(x-\frac{1}{x}\right) = \frac{-x^2 + 4x + 1}{x^2 + x - 1}=T(x)-\frac{1}{T(x)},$$
and we have the identity
\begin{equation}
(x+1)^5(y+1)^5 h(L(x),L(y)) = 5^5 h(y,x).
\label{eq:6.7}
\end{equation}
Moreover,
\begin{equation}
L(z)+1=\frac{5}{z+1}.
\label{eq:6.8}
\end{equation}

Using (\ref{eq:6.6}), (\ref{eq:6.7}) and (\ref{eq:6.8}), it follows by the same reasoning as in Section 2 that $S_n(x)$ has distinct roots and that
\begin{equation}
\textsf{Q}_n(x) = \prod_{k \mid n}{S_k(x)^{\mu(n/k)}}
\label{eq:6.9}
\end{equation}
is a polynomial.  Furthermore, all of the roots of $\textsf{Q}_n(x)$ lie in $\textsf{K}_5$.  From Theorem 6.3 we see that the polynomial $t_d(x)$ divides $\textsf{Q}_n(x)$ whenever the automorphism $\tilde{\tau}_5$ has order $n$, and from Theorem 6.2, we see that these are the only irreducible factors of $\textsf{Q}_n(x)$ over $\mathbb{Q}$.  This gives \bigskip

\noindent {\bf Theorem 6.5.} {\it For $n > 1$, the polynomial $\textsf{Q}_n(x)$ is given by the product
$$\textsf{Q}_n(x) = \pm \prod_{-d \in \mathfrak{D}_n^{(5)}}{t_d(x)},$$
where $t_d(x)$ is defined by (\ref{eq:6.5}) and $\mathfrak{D}_n^{(5)}$ is the set of negative quadratic discriminants $-d$ with $\left(\frac{-d}{5}\right)=+1$, for which the automorphism $\tilde{\tau}_{5,d} = \tilde{\tau}_5 = \left(\frac{\Omega_f/K}{\wp_5}\right)$ has order n in $\textrm{Gal}(\Omega_f/K)$, the Galois group of the ring class field $\Omega_f$ over $K=\mathbb{Q}(\sqrt{-d})$.} \medskip

For $\textsf{Q}_1(x)$ we have the factorization
\begin{align*}
\textsf{Q}_1(x) = & \ -(x + 1)(x^2 + 4)(x^2 - x + 3)(x^2 - 2x + 2)(x^2 + x + 5)\\
= & \ -(x+1)t_4(x) t_{11}(x) t_{16}(x) t_{19}(x),
\end{align*}
where $t_4(x)$ satisfies
$$x^2 t_4\left(x-\frac{1}{x}\right) = (x^2+1)^2=p_4(x)^2.$$
Since $\textrm{deg}(t_d(x)) = 2h(-d)$, Theorem 6.3 shows that half of the roots of $t_d(x)$ lie in the domain $\widetilde{\textsf{D}}_5$, while the other roots $\xi$ satisfy $\xi \equiv -1$ (mod $5$) in $\textsf{K}_5$, a fact which follows from (\ref{eq:6.7}) and (\ref{eq:6.8}).  Also see eq. (32) in \cite{mor8}. \medskip

The fact that $\textrm{deg}(t_d(x)) = 2h(-d)$ now implies the following class number formula.  \bigskip

\noindent {\bf Corollary 6.6.} {\it For $n > 1$ we have}
$$\sum_{-d \in \mathfrak{D}_n^{(5)}}{h(-d)} = \sum_{k \mid n}{\mu(n/k)5^k}.$$

This formula is equivalent to Deuring's formula for the prime $p=5$ from \cite{deu1}, \cite{deu3}, as in \cite{mor7}.

\medskip

\noindent Dept. of Mathematical Sciences, LD 270

\noindent Indiana University - Purdue University at Indianapolis (IUPUI)

\noindent Indianapolis, IN 46202

\noindent {\it e-mail: pmorton@iupui.edu}

\end{document}